\documentclass[12pt]{amsart}
\usepackage[utf8]{inputenc}
\usepackage[T1]{fontenc}

\usepackage{amsmath}
\usepackage{amssymb}
\usepackage{amscd}
\usepackage{amsfonts}
\usepackage{amsbsy}
\usepackage{graphicx}

\newcommand{\R}{\ensuremath{\mathbb{R}}}
\newcommand{\C}{\ensuremath{\mathbb{C}}}

\newcommand{\Z}{\ensuremath{\mathbb{Z}}}

\newtheorem {theorem} {Theorem} 
\newtheorem {proposition} [theorem] {Proposition}
\newtheorem {corollary} {Corollary}
\newtheorem {lemma} [theorem] {Lemma}
\newtheorem {definition} {Definition}

\begin{document}
\title[Perturbation theory of the Lotka-Volterra double center]
{Perturbation theory of the quadratic Lotka-Volterra double center}

\author[J.-P. Fran\c coise]
{Jean--Pierre Fran\c coise}

\address{Sorbonne-Universit\'e, Laboratoire Jacques--Louis Lions,
UMR 7598 CNRS, 4 Place Jussieu, 75252, Paris} 
\email{Jean-Pierre.Francoise@upmc.fr}

\author[L. Gavrilov]
{Lubomir Gavrilov}
\address{Institut de Math\'ematiques de Toulouse, Universit\'e de Toulouse, 31062, Toulouse, France}
\email{lubomir.gavrilov@math.univ-toulouse.fr}

\thanks{}


\dedicatory{}

\begin{abstract}
We revisit the bifurcation theory of the Lotka-Volterra quadratic system
\begin{eqnarray}
\label{X0}
X_0 :\left\{\begin{aligned}
\dot{x}=&   - y -x^2+y^2 ,\\
\dot{y}= &\;\;\;\;x - 2xy  
\end{aligned}
\right.
\end{eqnarray}
with respect to arbitrary quadratic deformations. The system (\ref{X0}) has a double center, which is moreover isochronous.  We show that the deformed system  (\ref{X0}) can have at most two limit cycles on the finite plane, with possible distribution $(i,j)$, where $i+j\leq2$. Our approach is based on the study of pairs of bifurcation functions associated to the centers, expressed in terms of iterated path integrals of length two.
\end{abstract}
\maketitle
\newpage
\tableofcontents
\section{Introduction}
\label{introduction}
We are interested in quadratic perturbations of  the following special reversible Lotka-Volterra quadratic system 
\begin{eqnarray}\label{nona}
X_0 :\left\{\begin{aligned}
\dot{x}=&   - y -x^2+y^2 ,\\
\dot{y}= &x - 2xy  
\end{aligned}
\right.
\end{eqnarray}
 \begin{figure}
    \includegraphics[height=8cm]{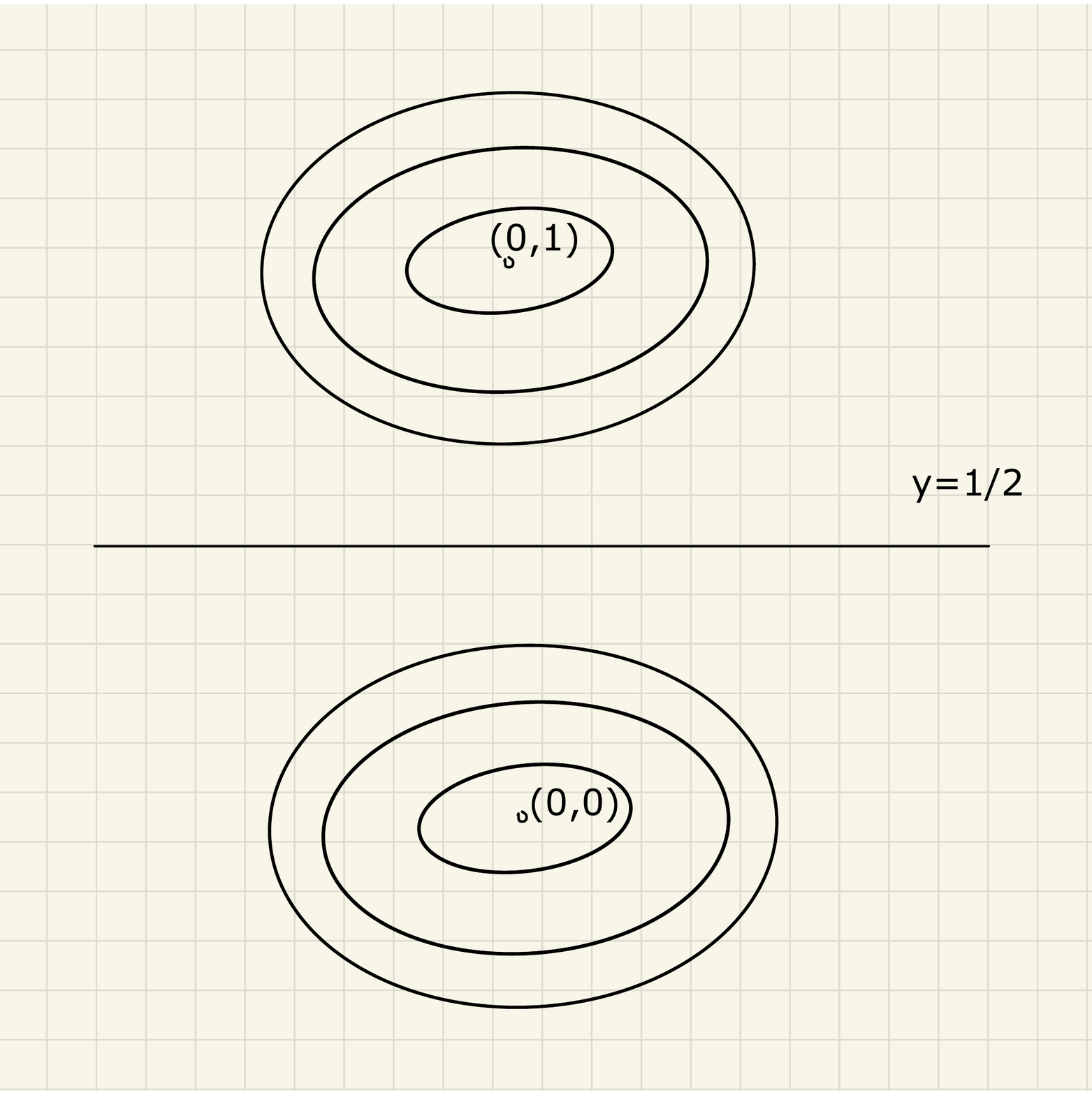} 
    \label{fig1n}
    \caption{Phase portrait of (\ref{nona}) }
     \end{figure}
which is equivalently written in  coordinates
 $z=x+iy$, $\bar z = x- i y$ as
\begin{equation}
\label{nonz}
\frac{dz}{dt} = iz-z^2 .
\end{equation}
This   implies that the vector field $X_0$ has a center at $z=0$ and $z=i$, that is to say at
the origin $(x,y)=(0,0)$ and at $(x,y)=(0,1)$, see  fig.1.
The period of the orbits is  
$$
T=\int dt = \oint \frac{dz}{ iz-z^2} = 2\pi 
$$
hence the two centers are
isochronous (the orbits have a constant period).

We are interested in the limit cycles which an arbitrary quadratic deformation of (\ref{nona}) can have. The limit cycles on a finite distance from the origin form two nests, containing either the focus close to $(0,0)$, or the focus close to $(0,1)$. We denote their number by $i$ and $j$.  The main result of the paper is easy to formulate: \emph{the possible distributions $(i,j)$ of limit cycles are those, for which $i+j \leq 2$.}

Although the above result is simple, it hides several difficulties, which were not resolved until recently. 
To the end of this Introduction, we outline the proof of Theorem \ref{main}, which will  also be
an occasion to illustrate some recent developments of the bifurcation theory of vector fields of infinite co-dimension.

The  system  (\ref{nona}) has, as suggested by (\ref{nonz}), a first integral
\begin{equation}
\label{hh}
H= \frac{x^2+y^2}{ 2y-1}  = \frac{x^2+(y-1)^2}{2y-1} +1
\end{equation}
 It induces a polynomial foliation on $\R^2$ (or $\C^2$) 
\begin{equation}
\label{nonperturbed2}
(1-2y)^2dH  = 0 .
\end{equation}
obviously invariant under the involution $x\mapsto -x$, but also with respect to the involution
\begin{align}
\label{involution}
y \mapsto 1-y .
\end{align}
The latter exchanges the two period annuli of the vector field $X_0$, which 
shows that the separate study of their deformations is analogous. If we prove that the cyclicity of the first period annulus is two, this implies that the cyclicity of the second period annulus is also two. The problem which we solve in the paper is the \emph{simultaneous} study of the bifurcations of the two period annuli.

An arbitrary quadratic perturbation of   (\ref{nonperturbed2}) or (\ref{nona}),  can be written in one of the following alternative forms
\begin{equation}
\label{pert1}
\frac 12 (1-y)^2dH + \sum_{0\leq i,j\leq 2 } (a_{ij} x^iy^j dy +  b_{ij} x^iy^j dx) = 0
\end{equation}
or
\begin{eqnarray}\label{pert2}
X_{a,b} :\left\{\begin{aligned}
\dot{x}=&   - y -x^2+y^2 + \sum_{0\leq i,j\leq 2 } a_{ij} x^iy^j ,\\
\dot{y}= & \;\;\;\;x - 2xy  -  \sum_{0\leq i,j\leq 2 }   b_{ij} x^iy^j 
\end{aligned}
\right.
\end{eqnarray}
In coordinates $z= x+ i y, \bar z = x- iy$, and
up to an affine transformation of $\R^2$ and a scaling of time, each of the above systems can be written in the following normal form
\begin{equation}
\label{perturbedz2}
\dot{z}=(\lambda_1+ i )z+ Az^2  + B z \bar z + C \bar z^2, B,C \in \C , \lambda_1\in \R
\end{equation}
where
\begin{align}
\label{ABC}
A=-1,\; B= \lambda_2+i \lambda_3, C= \lambda_4+i \lambda_5,\; \lambda_i \in \mathbb R.
\end{align}
and $\lambda_1,\lambda_2, \dots, \lambda_5$ are small real constants, see \cite{Z,I}.
We obtain finally the vector field
\begin{equation}
\label{perturbed2}
X_\lambda :\left\{
\begin{aligned}
\dot{x}=&  - y -x^2+y^2 
+ \lambda_1 x + \lambda_2 (x^2+y^2) + \lambda_4 (x^2-y^2) + 2\lambda_5 xy,\\
\dot{y}= &x- 2xy 
 + \lambda_1 y + \lambda_3 (x^2+y^2) + \lambda_5 (x^2-y^2) - 2 \lambda_4 xy.
\end{aligned}
\right.
\end{equation}
to be studied in the paper.
Thus, to obtain from (\ref{pert2}), the normal form (\ref{perturbed2}), we have to substitute

\begin{align*}
a_{10} &= \lambda_1,& b_{01}=&-\lambda_1\\
a_{20} &= \lambda_2+ \lambda_4, &b_{20} =& -\lambda_3-\lambda_5 \\
a_{02} &= \lambda_2-\lambda_4, & b_{02} = & -\lambda_3 + \lambda_5 \\
a_{11} &= 2 \lambda_5, & b_{11} =& 2 \lambda_4
\end{align*}
and $a_{00}=b_{00}=a_{01}=b_{10}=0$. 

The foliation underlying the vector field $X_\lambda$ takes the form
\begin{align}
\mathcal F_\lambda : \frac12 dH - \omega = 0
\end{align}
where
\begin{align}
\label{eq1}
\omega &= \frac{1}{(1-2y)^2} \sum_{0\leq i,j\leq 2 } (a_{ij} x^iy^j dy +  b_{ij} x^iy^j dx) = 
\sum_{i=1}^5 \lambda_i \omega_i&\\
\label{eq2}
\omega_1&=  \frac{xdy - ydx}{(2y-1)^2}, \omega_2 = \frac{(x^2  +y^2 )dy}{(2y-1)^2},
\omega_3 =-   \frac{ x^2+y^2}{(2y-1)^2} dx&    \\
\label{eq3}
\omega_4&= \frac{ (x^2-y^2)dy + 2xy dx}{(2y-1)^2} ,
\omega_5 = \frac{  2xy dy - (x^2-y^2)dx }{(2y-1)^2} .&
\end{align}

Let $\mathcal P(h)$ be the Poincaré first return map associated to one of the foci of (\ref{pert2}), which are close to $(0,0)$ and $(0,1)$ for $a_{ij}, b_{ij}$ sufficiently small. Here $h$ is as usual the restriction of the first integral $H$ of the non-perturbed system on  a transversal open segment through the focus. It is easily seen that $\mathcal P(h)$ is analytic both in $h$ and the parameters $a_{ij} , b_{ij} $, provided that the deformation is small and $h$ is close to the critical value of $H$. Expanding the displacement map $$\mathcal P(h)-h$$ in a power series in $h$ 
$$
\mathcal P(h)-h = \sum_{k=0}^\infty p_k(\lambda) h^k
$$
we consider the ideal 
$\mathcal B = < p_k(\lambda)> \subset  \R\{x,y\} $ generated by the coefficients of $h^k$. The fundamental fact about this Noetherian  ideal is, that it is \emph{ polynomially generated}.

The main advantage of the form (\ref{perturbedz2})  is that its Bautin ideal is known and relatively simple, which is not the case of (\ref{pert2}).
For this reason the  forms (\ref{perturbedz2}) and (\ref{perturbed2}) will be used from now on.

We denote $\mathcal B_1, \mathcal B_2 $ the local Bautin ideals associated to $(0,0)$ and $(0,1)$, localised at $\lambda=0$.
Our first result is the explicit form of the generators of $\mathcal B_1, \mathcal B_2$ in the paramer space $\R\{\lambda_1,\dots,\lambda_6\}$, see Theorem \ref{thbautin}. The proof uses in an essential way the explicit form of the so called Melnikov functions, combined with a version of the Nakayama lemma.
It follows from this result, that the irreducible algebraic set of quadratic systems of Lotka-Volterra type $\mathcal L(1,1,1)$, has a self-intersection at the "point" (\ref{X0}). This phenomenon is illustrated on Fig.\ref{self}. The two local branches of $\mathcal L(1,1,1)$ near the point (\ref{nonz}) are interchanged by the involution on the parameter space, induced by the affine involution $z \mapsto i-z$, see (\ref{involution}). As the cyclicity of $\mathcal L(1,1,1)$ is at least two, Fig.\ref{self} suggests the existence of a $(2,2)$ distribution, which was never met before in a quadratic system. This can be seen as \emph{a posteriori} motivation for the present paper. From this perspective our main result is negative. Indeed, neither $(2,2)$ nor even $(2,1)$ or $(1,2)$  distribution is possible on the finite plane under quadratic perturbation of $z'=iz-z^2$. 
\vskip 1pt
In a recent article, J.-P. Fran\c{c}oise, L. Gavrilov and D. Xiao (cf. \cite{FGX}) introduce an algebro-geometric setting for the space of bifurcation functions
involved in the local Hilbert’s 16th problem on a period annulus. Each possible bifurcation function is in one-to-one correspondence with a point in the exceptional divisor $E$ of the canonical blow-up $BI(\C^n)$ of the Bautin ideal $I$ . In this setting, the notion of essential perturbation,
first proposed by Iliev \cite{I}, is defined via irreducible components of the Nash space of arcs $Arc(BI(\C^n), E)$. The arcs here are seen as one-parameter deformations of a plane differential system with a center.
The example of planar quadratic vector fields with a center was further discussed. Here, in this article, we develop the same tools in the case of the quadratic Lotka-Volterra double center. In that case, the first order bifurcation function was recently directly computed by complex residues techniques in \cite{FY}. The first order Melnikov function had been computed by using the fact that (\ref{X0}) is isochronous with an explicit linearization in \cite{GGA} (see also \cite{LLLZ}). We mention previous important contributions of \cite{CJ} to the perturbation of Lotka-Volterra double center also based on Bautin ideal. The higher-order bifurcation function can, in general, be expressed in terms of iterated integrals (see \cite{F, G, FP}). We succeed here to compute the second-order bifurcation function using the shuffle formula and the complex residues techniques. An expression for this second-order bifurcation function had been found first by Iliev \cite[section 5 - III]{I}. This implies that the cyclicity of 
each annulus of $X_0$, see (\ref{X0}), is two, thus correcting a mistake in \cite{CJ}. 
Using the fact that the Lotka-Volterra double center is isochronous, C. Li and J. Llibre (\cite{LL2}) computed the second-order bifurcation function by averaging theory, by making use of   computer algebra. Partial results about the distribution of limit cycles were also obtained \cite[page 13]{LL2}. In particular, it was suggested that a $(2,1)$ distribution of limit cycles is possible, after suitable perturbation of $X_0$ on the finite plane. This claim is not correct. Note, however,  that nothing is known in this situation about the bifurcation of limit cycles from infinity.  We note also, that the cyclicity of the single period annulus of  a generic quadratic Lotka-Volterra system with three real invariant lines, is known to be two
by a classical result of {\.Z}o{\l}adek \cite{Z} (recently revisited in \cite{zola15}).

In contrast to the above papers,  we make essential use of the double Bautin ideal, see Theorem \ref{thbautin} and section \ref{sectionblowup}.
 We compute the first Melnikov functions related to the double ideal. This computation is possible, mainly because the complexified orbits $$\Gamma_h = \{ (x,y) \in \C^2 : H(x,y)=h\}$$ are genus zero curves. Thus, the first order Melnikov functions are computed by residue calculus, see sections \ref{firstmelnikov}, while for the second one we need shuffle relations of iterated integrals of length two, see section  \ref{secondmelnikov}. Alternatively, we describe in section \ref{section3} the geometry of the fibration defined by $H$, from which the monodromy of the bifurcation functions is explicitely determined. We determine in particular "the homology of the orbit", a new topological invariant adapted to the study of bifurcation functions. The corresponding section \ref{homology1} can be seen as an introduction to this notion and can be read independently. 
The description of the  monodromy of the bifurcation functions implies also their explicit form, section \ref{section4}, which provides an alternative proof of Theorem \ref{thm2}.
All this computations are performed for the center near $(0,0)$.

 The explicit expression that we obtain for both the two centers is very simple. It opens the way to compute completely the (double) Bautin ideal and prove that we can stop at order two. Hence, the final bound is deduced.

\vskip 1pt
For general references on Bifurcation Theory and Hilbert's 16th problem see \cite{CLW, R2}. After the classical contributions of NN Bautin \cite{B1,B2}, 
the notion of division by a Bautin ideal was first considered in \cite{FYo,R2}. References on the center-focus problem include \cite{CL,S}. The set of quadratic double centers was determined in \cite{L} and discussed more recently in \cite{FY}.

\section{The Bautin ideal and the bifurcation functions}
\label{bautinideal}
     \begin{figure}
    \includegraphics[height=6cm]{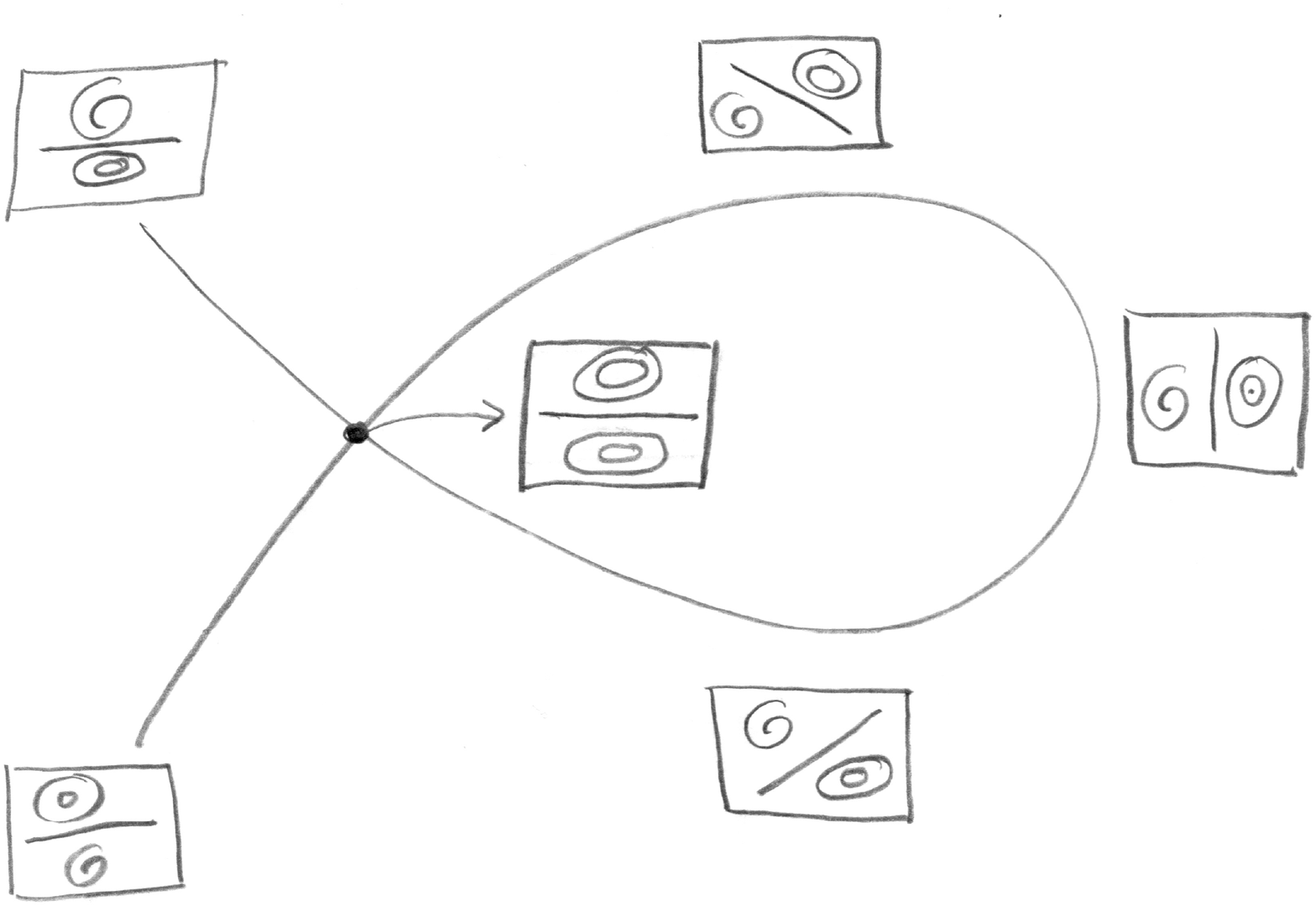} 
    \caption{The Lotka-Volterra component $\mathcal L (1,1,1)$ has a self-intersection at the isochronous system $\dot{z}=iz-z^2$. }
       \label{self}
     \end{figure}
       \begin{theorem}
     \label{thbautin}
 The vector field $X_\lambda$ (\ref{perturbed2}) has, for small parameters $\lambda_i$,  two foci close to $(0,0) $ and $(0,1)$.
The respective Bautin ideals $\mathcal B_1, \mathcal B_2$  are given by
\begin{align}
\label{bautin1}
\mathcal B_1 &= <\lambda_1,\lambda_3, \lambda_2\lambda_5 >\\
\label{bautin2}
 \mathcal B_2 &= <\lambda_1+\lambda_3 + \lambda_1 \lambda_2,\lambda_5, \lambda_3 \lambda_4>  .\end{align}
\end{theorem}    
  
    The zero locus of  $\mathcal B_1 $ has two irreducible components corresponding to systems or reversible or Lotka-Volterra type :
 \begin{align}
 \label{rv1}
 \lambda_1=\lambda_3= \lambda_5=0& \;\; \mbox{  (reversible component)}\\
 \label{lv1}
\lambda_1=\lambda_2= \lambda_3=0 &  \;\; \mbox{    (Lotka-Volterra component)}
 \end{align}
     with a similar structure of the isomorphic zero locus of $ \mathcal B_2$
     \begin{align}
 \label{rv2}
 \lambda_1=\lambda_3= \lambda_5&=0 \;\; \mbox{  (reversible component)}\\
 \label{lv2}
\lambda_1+ \lambda_3 + \lambda_1\lambda_2= \lambda_4=\lambda_5&=0
      \;\; \mbox{    (Lotka-Volterra component)}
   \end{align}
   Let us recall that the set of quadratic systems with a center form an algebraic subvariety (of the set of coefficients of all quadratic systems) with four irreducible components  \cite[Theorem 1.1.]{lins14} and \cite[Theorem 9]{FGX}. In particular, the Lotka-Volterra component, denoted further   $ \mathcal L(1,1,1)$ is an irreducible algebraic variety of co-dimension  three in the space of all quadratic vector fields.
   It follows from (\ref{lv1}) and (\ref{lv2}) that near $\lambda=0$ the component $ \mathcal L(1,1,1)$ has two local 
    irreducible components which intersect at the isochronous double center $\dot{z}=iz-z^2$ ($\lambda=0$).
     Therefore we obtain  
     \begin{proposition}
     The Lotka-Volterra irreducible component $\mathcal L (1,1,1)$ of the set of quadratic centers (the so called center set) has a self-intersection at the isochronous double center $\dot{z}=iz-z^2$, (\ref{nonz}).
     \end{proposition}
     The content of the above Proposition is illustrated on fig. \ref{self}. Note that in the Kapteyn-Dulac classification of quadratic centers, as presented by {\.Z}o{\l}adek \cite[Theorem 1]{Z}, the Lotka-Volterra stratum $Q_3^{LV}  $ is a plane, so it has no self-intersections. The reason is that  $Q_3^{LV}  $ is a quotient of $\mathcal L(1,1,1)$ with respect to the induced action of the affine group of linear changes of variables on $\R^2$.
Recall also that the general Lotka-Volterra system has a 
     single center (or Morse type equilibrium point) which is persistent, see \cite[Lins Neto ]{lins14}. The double centers of the isochronous system
        $\dot{z}=iz-z^2$ are in this sense semi-persistent. This means that each of the centers is persistent or non-persistent depending on the local branch of the stratum  $\mathcal L (1,1,1)$. This possibility is not mentioned in  \cite{lins14}.

\begin{proof}[Proof of Theorem \ref{thbautin}]
 The proof of (\ref{bautin1}) goes back to Dulac (1908) \cite{D}.
     Indeed, following \cite[{\.Z}o{\l}adek]{Z} and \cite[I, pp.122-123]{I} we deduce that $\mathcal B_1$ is generated by 
     $\lambda_1$ and the focal values $v_3,v_5,v_7$ where
     \begin{align}
     v_3 &= 2\pi \mathrm{Im} AB  \\
     v_5 &=  \frac23 \mathrm{Im} [(2A+\bar B)(A-2\bar B) \bar B C] \\
     v_7 &= \frac54 (|B|^2-|C|^2) \mathrm{Im} [(2A+ \bar B) \bar B^2C] .
     \end{align}
 According to (\ref{ABC})    $A=-1$ and $v_3=-2\pi \lambda_3$. Assuming that $\lambda_3 = 0$ we have
      \begin{align}
      v_5 & =  \frac23 \mathrm{Im} [(-2+\lambda_2)(-1-2\lambda_2) \lambda_2 ( \lambda_4+i \lambda_5)] \\
      &= \frac23 (-2+\lambda_2)(-1-2\lambda_2) \lambda_2  \lambda_5
\end{align}
     and therefore locally $<v_3,v_5> = < \lambda_3, \lambda_2 \lambda_5>$ and $v_7$ is generated by $v_3,v_5$, which proves (\ref{bautin1}).
 
The proof of (\ref{bautin2}) will be done in two steps. 
\begin{enumerate}
\item 
First, we prove that the center set of the second center (the variety of the ideal $\mathcal B_2$) is defined by
(\ref{rv2}), (\ref{lv2}). For this, we show that when the vector field satisfies (\ref{rv2}), (\ref{lv2}), then it has a first integral analytic near $(0,1)$ and hence has a center. 
\item
Second, we show that the   ideal $\mathcal B_2$ is radical. For this we use the information obtained from the computation of the first and the second Melnikov functions, and a version of Nakayama lemma, as suggested by \cite[Lemma 7.4]{bry10}.
\end{enumerate}
The involution $z\to -z+i$  exchanges the singular points near $(0,0)$ and $(0,1)$ and acts linearly on the coefficients of the vector field $X_{a,b}$. Therefore the two center sets are analytically isomorphic (near $\lambda=0$). It follows that the germ of variety of $\mathcal B_2$ has two smooth irreducible components of co-dimension three. The first one is obviously the reversible one (\ref{rv2}), and the system has two reversible double centers. To show that (\ref{lv2}) is the second irreducible component, we shall find 
a  first integral of $X_\lambda$, which is analytic near the point $(0,1)$.

In what follows we assume that $\lambda_5 = \lambda_4 = 0 $.
The equation (\ref{perturbedz2}) takes the form
\begin{align}
\label{second}
\dot{z}=(\lambda_1+ i )z-z^2  + B z \bar z = z( \lambda_1+ i +B \bar z)
\end{align}
where $B=\lambda_2+i \lambda_3 \in \C $ and  $ \lambda_1, \lambda_2, \lambda_3 \in \R$, and
  can be  integrated as follows. Consider  the underlying foliation defined by
\begin{align}
\label{fsecond}
(1- i \lambda_1 + i z - iB \bar z) z d \bar z + (1+ i \lambda_1 - i \bar z + i\bar B  z) \bar z d  z = 0 .
\end{align}
If it were integrable, it would be of Lotka-Volterra type, and therefore
would have at least three invariant lines intersecting at singular points of the foliation. Indeed, the following two lines are obviously invariant
\begin{align*}
z=x+iy=0, \bar z = x- iy=0
\end{align*}
and let the third one be
\begin{align*}
\alpha z + \bar \alpha \bar z +1 = 0, \alpha \in \C .
\end{align*}
This implies on its turn an ansatz for the  first integral as follows :
\begin{align*}
H = z^{1+ i \lambda_1}\bar z ^{1- i \lambda_1}(\alpha z + \bar \alpha \bar z +1 )^\beta, \beta \in \R .
\end{align*}
The polynomial foliation $ d \log H  = 0$ is
\begin{align*}
(1+ i \lambda_1) \frac{dz}{z} + (1- i \lambda_1) \frac{d\bar z}{\bar z} + \beta \frac{d ( \alpha z + \bar \alpha \bar z +1)}{\alpha z + \bar \alpha \bar z +1}=0
\end{align*}
or equivalently
\begin{align*}
( \alpha z + \bar \alpha \bar z +1) [(1+ i \lambda_1)\bar z  dz + (1- i \lambda_1) z  d\bar z ]
+\beta z\bar z d ( \alpha z + \bar \alpha \bar z +1) &= 0 \\
(1+ i \lambda_1)\bar z  dz + (1- i \lambda_1) z  d\bar z  
+ (1+ i \lambda_1)(\alpha z + \bar \alpha \bar z) \bar z  d  \\
+ (1- i \lambda_1)( \alpha z + \bar \alpha \bar z) z  d\bar z 
+  \beta z \bar z d ( \alpha z + \bar \alpha \bar z ) &= 0
\end{align*}
and finally
\begin{align*}
(1- i \lambda_1) z  d\bar z + [(1- i \lambda_1) \alpha z^2 + (\bar \alpha (1- i \lambda_1) + \beta \bar \alpha )z\bar z ] d \bar z & \\
+ (1+ i \lambda_1) \bar z  d z + [(1+ i \lambda_1) \alpha \bar z^2 + ( \alpha (1+ i \lambda_1) + \beta  \alpha )\bar z z ] d  z &=0 .
\end{align*}
Comparing this to (\ref{fsecond}) we impose
\begin{align*}
(1- i \lambda_1) \alpha & = i \\
\bar \alpha (1- i \lambda_1 + \beta   ) &= - i B
\end{align*}
where $B= \lambda_2+i \lambda_3$. Therefore
\begin{align*}
1- i \lambda_1 + \beta  = (\lambda_2+i \lambda_3) (1+ i \lambda_1)
\end{align*}
and finally
\begin{align*}
1+ \beta &= \lambda_2 -\lambda_1 \lambda_3 \\
-\lambda_1 &= \lambda_1 \lambda_2 + \lambda_3 .
\end{align*}
The conclusion is that if 
$$
\lambda_1 + \lambda_3 + \lambda_1 \lambda_2  = 0
$$
then
$$
H = \frac{z^{1+ i \lambda_1}\bar z ^{1- i \lambda_1}}{(\alpha z + \bar \alpha \bar z +1 )^{1-\lambda_2 + \lambda_1 \lambda_3 }},
\alpha = \frac{i-\lambda_1}{1+\lambda_1^2}
$$
is a first integral of (\ref{second}). In other words, the variety 
\begin{align*}
\{\lambda\in \C^5 : \lambda_1 + \lambda_3 + \lambda_1 \lambda_2  = \lambda_4=\lambda_5=0 \}
\end{align*}
is a co-dimension   three irreducible component of the center manifold related to the second center. This completes the first step of our proof.

In the ring of convergent power series $\R\{\lambda\}$ consider the ideal of functions vanishing along the variety
(\ref{rv2}), (\ref{lv2}). It is obviously generated by
$$
a= \lambda_1+\lambda_3 + \lambda_1 \lambda_2, b= \lambda_5, c=  \lambda_3 \lambda_4 
$$
and at the second step we shall show that
$$
\mathcal B_2 = <a,b,c>.
$$
We examine first the information obtained from the first and the second Melnikov functions (see the next sections).
It follows from (\ref{m1t}) that there are elements  $ v_1^2,  v_2^2$ of the ideal $\mathcal B_2$ such that
\begin{align*}
 v_1^2(\lambda)& = \lambda_1+\lambda_3  + \dots \\
 v_2^2(\lambda) &= \lambda_5 + \dots
\end{align*}
where the dots replace some analytic series which vanish of order at least two at $\lambda=0$. 
We can write therefore
\begin{align*}
 v_1^2(\lambda)& = a + \alpha c +\dots \\
 v_2^2(\lambda) &= b + \beta c + \dots
\end{align*}
where the dots replace  some analytic series which vanish along (\ref{rv2}), (\ref{lv2}),
 vanish of order at least two at $\lambda=0$, and 
$\alpha,\beta $ are appropriate constants.
We can write finally
\begin{align}
\label{v1}
 v_1^2(\lambda)& = a\,(1+O(\lambda )) + b \,O(\lambda)   + c \,( \alpha + O(\lambda)) \\
 \label{v2}
 v_2^2(\lambda) &= a \,O(\lambda) + b\,(1+ O(\lambda))  + c \,(\beta + O(\lambda) ) .
\end{align}
Similarly, the identity (\ref{97}) implies that under the condition
$$ \lambda_1+\lambda_3 =   \lambda_5 = 0$$ there is an element $v_3^2$ of $\mathcal B_2$ such that
\begin{align*}
 v_3^2(\lambda)& =  \lambda_3\lambda_4  + \dots 
\end{align*}
where the dots replace some analytic series vanishing along (\ref{rv2}), (\ref{lv2}), and 
vanish of order at least three at $\lambda=0$.
Without the conditions $ \lambda_1+\lambda_3 =  \lambda_5 = 0$, we get
\begin{align*}
 v_3^2(\lambda)& =  \lambda_3\lambda_4  + ( \lambda_1+\lambda_3) (\gamma + O(\lambda)) +
 \lambda_5 (\delta+ O(\lambda)) + ...
\end{align*}
where the dots replace some analytic series which vanish of order at least three at $\lambda=0$.
Thus
\begin{align*}
 v_3^2(\lambda)& = c (1+ O(\lambda) ) + a\, (\gamma + O(\lambda)) +
b \,(\delta+ O(\lambda)) .
\end{align*}
and combining with (\ref{v1}), (\ref{v2})
\begin{align}
\left(\begin{array}{c} v_1^2\\v_2^2 \\v_3^2\end{array}\right)
=
\left(\begin{array}{ccc}1+ O(\lambda) & 0 & \alpha+ O(\lambda) \\0 & 1+ O(\lambda) & \beta+ O(\lambda) \\
\gamma + O(\lambda)& \delta +  O(\lambda) & 1+  O(\lambda)\end{array}\right) \left(\begin{array}{c} a\\b \\c\end{array}\right)
\end{align}
As the above matrix is invertible for $\lambda$ close to the origin, then $a,b,c$ belong to $\mathcal B_2$, which
 completes the proof of  Theorem \ref{thbautin} .
\end{proof}

The vector field  $X_\lambda$, see (\ref{perturbed2}), defines  return maps $\mathcal P_1, \mathcal P_2$ with associated Bautin ideals
\begin{align}
\label{bb1}
\mathcal B_1 &= < v_1^1(\lambda),v_2^1(\lambda),v_3^1(\lambda)> =<\lambda_1,\lambda_3, \lambda_2\lambda_5 >\\
\label{bb2}
\mathcal B_2 &=  < v_1^2(\lambda),v_2^2(\lambda),v_3^2(\lambda)> =<\lambda_1+\lambda_3 + \lambda_1 \lambda_2,\lambda_5, \lambda_3 \lambda_4> .
\end{align}
 $\mathcal P_1, \mathcal P_2$  
 can be divided 
in the corresponding ideals (\ref{bb1}) and (\ref{bb2}) as follows, (see \cite{FYo}, \cite[section 4.3]{R}).

\begin{align*}
\mathcal P_1^1(h;\lambda)(h) - h &=  v_1^1(\lambda)(M_1^1(h) + O(\lambda)) + v_2^1(\lambda)(M_2^1(h) 
+ O(\lambda)) \\ 
&+ v_3^1(\lambda)(M_3^1(h) + O(\lambda))\\
\mathcal P_1^2(h;\lambda)(h) - h &= v_1^2(\lambda)(M_1^2(h) + O(\lambda)) + v_2^2(\lambda)(M_2^2(h) + O(\lambda)) \\
&+v_3^2(\lambda)(M_3^2(h) + O(\lambda)). 
\end{align*}
\begin{definition}
\label{bf1}
The functions 
$$
M_1^1,  M_2^1 \mbox{ and } M_1^2, M_2^2
$$
are called the first order (or linear) Melnikov functions, associated to the centers at $(0,0)$ and $(0,1)$.
The functions 
$$
M_3^1 \mbox{ and } M_3^2
$$
are called the second order (or non-linear) Melnikov functions, associated to the centers at $(0,0)$ and $(0,1)$.
\end{definition}
The terminology is due to  \cite[{\.Z}o{\l}adek, section 2]{Z} and it
will be justified in what follows.
Given an arc,
\begin{align}
\label{arc}
\varepsilon \mapsto \lambda(\varepsilon), \; \varepsilon \in (\R,0), \;\lambda(0)=0, 
\end{align}
 we obtain

\begin{align*}
\mathcal P_1^1(h;\lambda(\varepsilon))(h) - h &= \varepsilon^{k_1}( c_1^1 M_1^1(h)  + c_2^1 M_2^1(h) 
+   c_3^1 M_3^1(h) + O(\varepsilon))\\
\mathcal P_1^2(h;\lambda(\varepsilon))(h) - h &=  \varepsilon^{k_2}( c_1^2 M_1^2(h)  + c_2^2 M_2^2(h) 
+   c_3^2 M_3^2(h) + O(\varepsilon)).
\end{align*}
Note that not all linear combinations 
\begin{align}
\label{admissible}
c_1^1 M_1^1(h)  + c_2^1 M_2^1(h) 
+   c_3^1 M_3^1(h), c_1^2 M_1^2(h)  + c_2^2 M_2^2(h) 
+   c_3^2 M_3^2(h)
\end{align}
of Melnikov functions are admissible. 
\begin{definition}
Let $K\subset \R^2$ be a compact set. A $(i,j)$ distribution  of limit cycles is said to be admissible for $X_\lambda$, if for every $\varepsilon >0$ there exists $\lambda$, such that
$\|\lambda \| < \varepsilon$ and $X_\lambda$ has a  $(i,j)$ distribution of limit cycles in $K$. 
\end{definition} 
Let  $(i,j)$ be admissible distribution of limit cycles for
$X_\lambda$ in the compact set $K$. Then there exists a germ of analytic arc
(\ref{arc}),
such that the one-parameter family of vector fields $X_{\lambda(\varepsilon)}$ allows a distribution $(i,j)$  limit cycles, for  $\varepsilon$ close to $0$,  \cite[Theorem 1]{G08}. Therefore to compute the possible distributions $(i,j)$ of limit cycles we have to compute the number of zeros $i$ and $j$
of each admissible pair of Melnikov functions (\ref{admissible}).

Denote
\begin{align*}
\delta(h) &= \{(x,y)\in \R^2 : H(x,y) = h\} , h \leq 0 \\
\tilde \delta(h) &= \{(x,y)\in \R^2 : H(x,y) = h\} , h \geq 1 
\end{align*} 
the family of real ovals  of the affine algebraic curve $\Gamma_h \in \C^2$. Using the notations (\ref{eq1}), (\ref{eq2}), (\ref{eq3}) ,
the following integral formulae for the linear Melnikov functions are well known,

\begin{proposition}
The linear Melnikov functions are given by
\begin{align*} 
\frac12 M^1_1(h) =  \int_{\delta(h)}   \omega_1, &\quad \frac12 M^1_2(h) =  \int_{\delta(h)}   \omega_3 \\
\frac12 M_1^2(h)  =  \int_{\tilde \delta(h)}   \omega_1, &\quad  \frac12 M_2^2(h)  =  \int_{\tilde \delta(h)}   \omega_5 .
\end{align*}
\end{proposition}
\proof
It is easy to verify that
\begin{align*} 
 \int_{\delta(h)}   \omega_2 &=  \int_{\delta(h)}   \omega_4 =  \int_{\delta(h)}   \omega_5 = 0 \\
 \int_{\tilde \delta(h)} \omega_2 &= \int_{\tilde \delta(h)} \omega_4 = 0, \int_{\tilde \delta(h)}   \omega_1=\int_{\tilde \delta(h)}   \omega_3 .
\end{align*} 
and hence
\begin{align*} 
\int_{\delta(h)}   \omega = \lambda_1 \int_{\delta(h)}   \omega_1 + \lambda_3 \int_{\delta(h)}   \omega_3 \\
\int_{\tilde \delta(h)}   \omega = (\lambda_1+\lambda_3) \int_{\tilde \delta(h)}   \omega_1  +  \lambda_5 \int_{\tilde \delta(h)}   \omega_5 .
\end{align*} 
\endproof
The second order (nonlinear) Melnikov function is given by the second order in $\lambda$ homogeneous piece of the displacement maps  $\mathcal P_1-id, \mathcal P_2 - id$.
For a differential one-form on $\C^2$ let $\omega'$ be the Gelfand-Leray residue of $\omega$ with respect to $H$ defined by the identity
$$
\omega' \wedge dH = d \omega .
$$
The second order Melnikov function of a deformed foliation $dH+\varepsilon \omega=0$ is defined by the following iterated integral of length two \cite{G}
$$
\int_{\delta(h)}  \omega\omega' 
$$
with appropriate choice of the path $\delta(h)$. In our case this implies
\begin{proposition}
Assume that the linear Melnikov function $M^1_1=M^1_2=0$ . Then
\begin{align}
\frac14 M_3^1(h) = \int_{\delta(h)} \omega_2 \omega_5' + \omega_5 \omega_2'
\end{align}
Similarly, if $M^2_1=M^2_2=0$, then
\begin{align}
\frac14 M_3^2(h) = \int_{\tilde\delta(h)} (\omega_3-\omega_1) \omega_5' + \omega_5 (\omega_3-\omega_1)'
\end{align}\end{proposition}
\proof
The vanishing of $ \int_{\delta(h)}   \omega$ implies $\lambda_1=\lambda_3 = 0$.
According to \cite{G} the function $M_3^1(h)$ corresponds to the coefficient  $\lambda_2\lambda_5$ in the iterated integral
$$
\int_{\delta(h)} \omega \omega' , \omega = \sum_{i=1}^5 \lambda_i \omega_i .
$$
Therefore, assuming in addition that $\lambda_4=0$ we get

\begin{align*}
\int_{\delta(h)} \omega \omega' &= \int_{\delta(h)} (\lambda_2\omega_2+\lambda_5 \omega_5)(\lambda_2\omega_2+\lambda_5 \omega_5)'\\
& = \lambda_2\lambda_5 \int_{\delta(h)}  \omega_2 \omega_5' + \omega_5 \omega_2' 
\end{align*}
where we used that $ \int_{\delta(h)} \omega_2 \omega_2'=  \int_{\delta(h)} \omega_5 \omega_5'=0$ (This will be justified latter in the text by using the shuffle formula).
The proof of the formula for $M_3^2(h)$ follows the same lines.
\endproof

\section{The first Melnikov function}
\label{firstmelnikov}
In this section we compute, for completeness the first Melnikov functions of $X_\lambda$. These results are classical, see  \cite{CJ,GGA,I,LLLZ, LL2}.
Here we use a simple residue calculus, following  \cite{FY}. If we write  $X_\lambda$ in the form (\ref{pert1})

\begin{align}
\frac 12 (1-2y)^2dH + \omega = 0
\end{align}
where
$$
H= \frac{x^2+y^2}{2y-1}
$$

  Denote
\begin{equation}
\label{gammah}
 \Gamma_h = \{(x,y) \in \C^2 :  x^2+y^2 = (2y-1)h, y\neq \frac12 \} .
\end{equation}
which, for $h\neq 0,1$, 
is a four-punctured Riemann sphere, where the punctures are at 
\begin{align}
\label{punctures}
(\pm \frac{\sqrt{-1}}{2}, \frac 12 ) , \infty^\pm .
\end{align}
Let 
\begin{align}
\delta(h), \tilde\delta(h) \in H_1(\Gamma_h, \Z)
\end{align}
 be a continuous family of cycles vanishing at the singular points $(0,0)$ and $(0,1)$, when $h$ tends to $h=0$ or $h=1$ respectively. These two families of cycles are defined in a neighbourhood of $h=0$ and $h=1$ respectively, and hence on the real segment $(0,1)$.

\begin{definition}
The first Melnikov functions $M_1$,  $\tilde M_1$ associated to the centers $(0,0)$ and $(0,1)$ respectively,
are defined by 
$$M_1(h)= 2\int_{\delta(h)} \omega, \; \tilde M_1(h) =  2\int_{\tilde\delta(h)} \omega . $$ 
\end{definition}
The functions are analytic on $(0,1)$ and therefore can be computed and compared there. This is easy, as they are Abelian integrals on a Riemann sphere, so the computation is reduced to residue calculus. 

\begin{figure}
    \includegraphics[height=6cm]{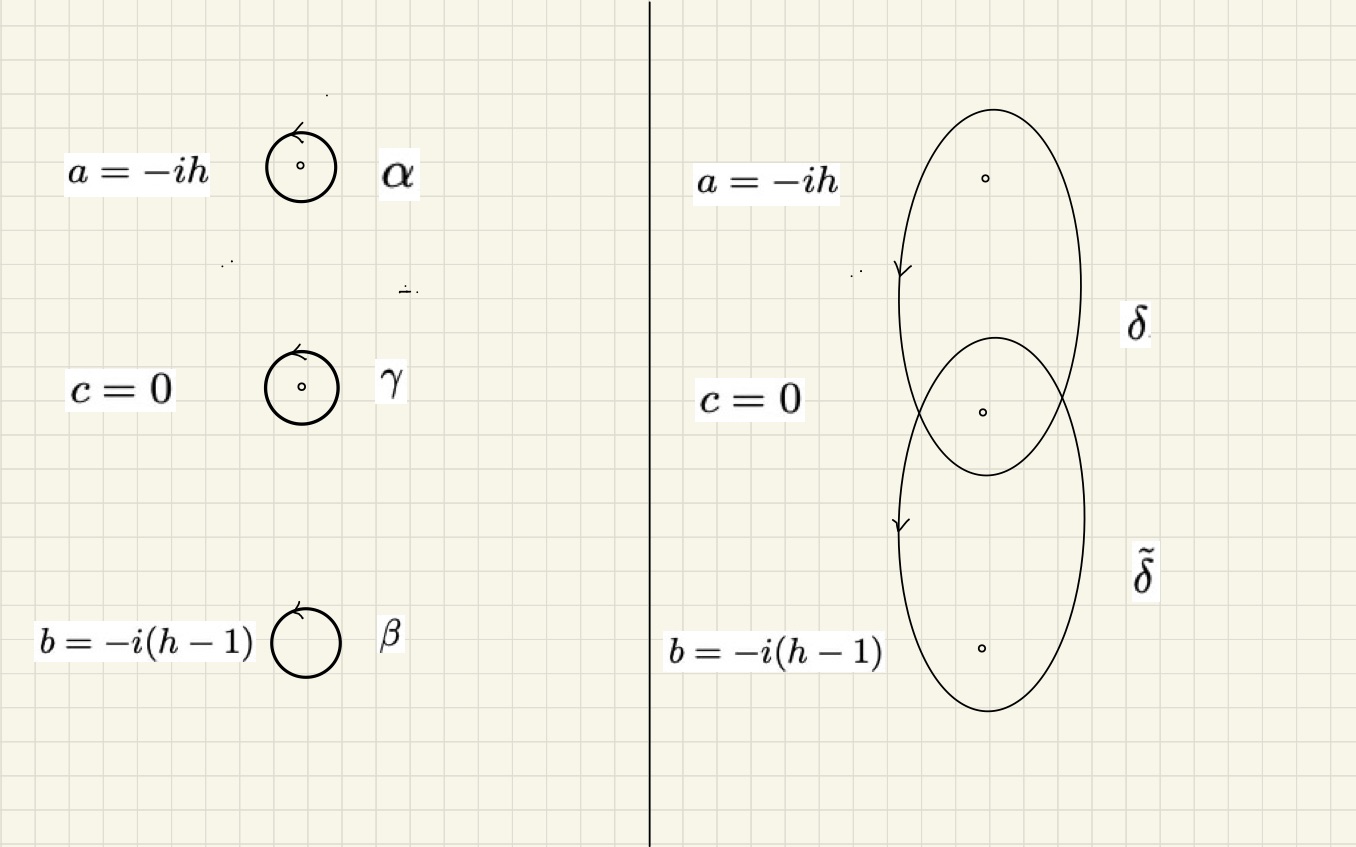} 
    \caption{The loops $\alpha, \beta, \gamma, \delta$ and $\tilde \delta$ for $h\in (0,1)$. }
       \label{figzf}
     \end{figure}

Following  \cite{FY}, chose an uniformizing variable $z: \bar \Gamma_h \to \mathbb P^1$ by the formula
\begin{align}
z= x + i(y-h), i= \sqrt{-1} .
\end{align}
If we note $\bar z = x - i(y-h)$ (so that $\bar z$ is complex conjugate to $z$ when $h\in \R$) we have
$$
\Gamma_h = \{(z,\bar z) \in \C^2 : z \bar z = h(h-1) \}.
$$
The images of the  four punctures  (\ref{punctures}) on the curve (\ref{gammah}) under $$z: \bar \Gamma_h \to \mathbb P^1$$
are 
\begin{align}
z(\infty^+) = \infty, z(\infty^-) = 0, z(\frac i2, \frac 12)= - i(h-1), z(-\frac i2, \frac 12)= - i h  
\end{align}
where $i$ is an appropriate determination of $\sqrt{-1}$. The model of the four-punctured Riemann sphere $\Gamma_h$ will be therefore the punctured complex plane
$\C \setminus \{a,b,c\}$, where
$$
a= -ih, b =  -i(h-1), c = 0 .
$$
$$
\alpha, \beta , \gamma, \delta, \tilde \delta
$$
We have
$$
\lim_{h\mapsto 0} a(h)= c, \; \lim_{h\mapsto 1} b(h)= c
$$
and it is easy to check that the vanishing cycles $\delta(h)$ and $\tilde \delta(h)$ are represented by "small" simple loops containing $a,c$ for $h\sim 0$, and $b,c$ for $h\sim 1$, as shown on fig.\ref{figzf}. It follows that for the homology classes (denoted by the same letters) holds 
$$
\delta = \alpha + \gamma, \; \tilde \delta = \beta + \gamma 
$$
and hence
$$
\frac12 M_1(h)= \int_{\alpha(h)} \omega + \int_{\gamma(h)} \omega, \; \frac12 \tilde M_1(h)= \int_{\beta(h)} \omega + \int_{\gamma(h)} \omega .
$$
The explicit computation of $M_1$ is a simple residue calculus. It was already computed by Fran\c{c}oise and Yang \cite{FY} and we reproduce them below.

With the notations (\ref{hh}), (\ref{eq1}), (\ref{eq2}), (\ref{eq3})
it  follows from \cite{FY} that
\begin{align*}
M_1(h) &=  -\frac{1}{16} \pi 4h(A_1\frac{4h-1}{8}+A_0), & h< 0 \\
\tilde M_1(h) &=  \frac{1}{16} \pi 4(h-1)(B_1\frac{4h-1}{8}+B_0), & h> 1 \\
\end{align*}
where 
\begin{align*}
A_1 &= 16(\lambda_3+\lambda_1),&           A_0 &= 2(\lambda_3-3\lambda_1) \\
B_1 &= A_1, &                             B_0 &= 2(\lambda_3+\lambda_1) - 16 \lambda_5
\end{align*}
which implies
\begin{align}
\label{m1}
M_1(h) &=&  -2\pi h [h(\lambda_1+\lambda_3) -\lambda_1] &, & h< 0 \\
\label{m1t}
\tilde M_1(h) &=&  2\pi (h-1) [(h-1)(\lambda_1+\lambda_3)+\lambda_1+\lambda_3-2\lambda_5] &, & h> 1 .
\end{align}
As expected  $M_1 \neq \tilde M_1$ which allows to construct at a first order all possible distributions $(i,j)$ of limit cycles, such that $i\leq 1, j\leq 1$.

\section{The second Melnikov function : analytic computation}
\label{secondmelnikov}
Let
$$H=\frac{x^2+y^2}{(2y-1)}$$
be the first integral of $X_0$, see (\ref{nona}).
On each level set 
$$\Gamma_h = \{(x,y)\in \C^2: H(x,y) =h\}
$$ holds
$$x^2+y^2=h(2y-1),$$
$$x^2+(y-h)^2=h(h-1).$$
The real level sets $\{(x,y)\in \R^2 : H(x,y)=h\}$, $h\in \R$, are therefore 
 circles centered at $(0,h)$ of radius $R=\sqrt{h(h-1)}$. The critical values of $H$ are $h=0$ and $h=1$. The two period annuli are 
 $$\{(x,y)\in \R^2 : H(x,y)<0\}, \;\;\{(x,y)\in \R^2  : H(x,y)>1\}.$$ Note that the symmetry $\sigma: y\to 1-y$ induces $\sigma^*(H)=H-1$. 
 Recall from section \ref{firstmelnikov} that $\Gamma_h$ is a four-punctured Riemann sphere, uniformized by
the complex parameter 
$$z = x+{\rm i}(y-h)$$
where
$$x-{\rm i}(y-h)=R^2/z.$$
$\Gamma_h$ is therefore identified with the complex $z$-plane with three punctures at $a, b, c$ where
$$
a=-{\rm i}h, b=-{\rm i}(h-1), c= 0 \mbox{  and  } R^2=-ab 
$$
see Fig. \ref{figzf}.
In what follows, as in the preceding section, $\omega$ is the differential one-form (\ref{eq1}), but under the condition that
$$
\int_{\delta(h)} \omega = 0 
$$
or
$$
\int_{\tilde \delta(h)} \omega = 0 .
$$
The loops $\delta, \tilde \delta$ are represented by circles surrounding $a, c$ or $b, c$ respectively, see Fig.\ref{figzf}. The one-form $\omega$ is holomorphic on $\Gamma_h$ and has poles at $z= a, b, c, \infty$.

Our purpose is to compute the second Melnikov function of the perturbed equation $dH-\omega = 0$ (the $1/2$ factor of $H$ was skipped for convenience).

\subsection{Computation of the perturbative part  in new coordinates}
In the normal form, the perturbative part can be written as:
\begin{equation}
\begin{array}{l}
\omega=\lambda_1\frac{xdy-ydx}{(2y-1)^2}+\lambda_2\frac{(x^2+y^2)dy}{(2y-1)^2}-\lambda_3\frac{(x^2+y^2)dx}{(2y-1)^2}\\
+\lambda_4\frac{(x^2-y^2)dy+2xydx}{(2y-1)^2}+\lambda_5\frac{2xydy-(x^2-y^2)dx}{(2y-1)^2}.
\end{array}
\end{equation}
If we assume that $M_1(h)\equiv 0$, then $\lambda_1=\lambda_3=0$ so
\begin{align*}
\omega &=\lambda_2\omega_2+\lambda_4\omega_4+\lambda_5\omega_5 \\
&=\lambda_2\frac{(x^2+y^2)dy}{(2y-1)^2}+\lambda_4\frac{(x^2-y^2)dy+2xydx}{(2y-1)^2}\\
&+\lambda_5\frac{2xydy-(x^2-y^2)dx}{(2y-1)^2}.
\end{align*}
It is easily verified that 
\begin{align*}
\C^2 \to \C^2 : (x,y) \mapsto (z,h), \;\; h=H(x,y)) 
\end{align*}
is a bi-rational transformation of $\C^2$.
Therefore we can use $z,h$ coordinates to express $dH-\omega$ and compute the corresponding second Melnikov function.
We have
$$dx=\frac{(z^2-R^2)}{2z^2}dz+\frac{1}{2}\frac{2h-1}{z}dh$$
$$dy=\frac{1}{2{\rm i}z^2}(z^2+R^2)dz+\frac{[z-\frac{a+b}{2}]}{z}dh.$$
We get by simple substitutions:
\begin{align*}
\frac{x^2dy}{(2y-1)^2}&=\frac{{\rm i}}{8}\frac{(z^2-ab)^3}{[z(z-a)(z-b)]^2}dz\\
&-\frac{1}{4}\frac{(z^2-ab)^2(z-\frac{a+b}{2})}{z[(z-a)(z-b)]^2}dh.
\end{align*}

\begin{align*}
\frac{y^2dy}{(2y-1)^2} =&-\frac{{\rm i}}{8}\frac{[(z-a)(z-b)
+{\rm i}z]^2(z^2-ab)}{[z(z-a)(z-b)]^2}dz\\
&+\frac{1}{4}\frac{[(z-a)(z-b)+{\rm i}z]^2(z-\frac{a+b}{2})}{z[(z-a)(z-b)]^2}dh.
\end{align*}

\begin{align*}
\frac{2xy}{(2y-1)^2}dx=&-\frac{1}{4{\rm i}}\frac{(z^2-ab)(z^2+ab)[(z-a)(z-b)+{\rm i}z]}{[z(z-a)(z-b)]^2}dz \\
& -\frac{(2h-1)}{4{\rm i}}\frac{(z^2-ab)[(z-a)(z-b)+{\rm i}z]}{z[(z-a)(z-b)^2}dh.
\end{align*}
\begin{align*}
\frac{2xy}{(2y-1)^2}dy=&\frac{1}{4}\frac{(z^2-ab)^2[(z-a)(z-b)+{\rm i}z]}{[z(z-a)(z-b)]^2}dz\\
&-\frac{1}{2{\rm i}}\frac{(z^2-ab)(z-\frac{a+b}{2})[(z-a)(z-b)+{\rm i}z]}{z[(z-a)(z-b)]^2}dh.
\end{align*}
\begin{align*}
\frac{-x^2+y^2}{(2y-1)^2}dx=&\frac{1}{8}\frac{((z^2-ab)^2+[(z-a)(z-b)+{\rm i}z]^2)(z^2+ab)}{[z(z-a)(z-b)]^2}dz\\
&+\frac{2h-1}{8}\frac{(z^2-ab)^2+[(z-a)(z-b)+{\rm i}z]^2}{z[(z-a)(z-b)]^2}dh.
\end{align*}
\subsection{The second-order Melnikov function defined by an iterated integral}
From $\omega=Fdz+\Phi dH,$ we get:
$$d\omega=(F'_H-\Phi'_z)dH{\wedge}dz.$$
The Gelfand-Leray derivative of $\omega$ is defined (modulo $dH$) by
\begin{equation*}
\omega'=(F'_H-\Phi'_z)dz.
\end{equation*}
The associated second-order Melnikov function is defined as the iterated integral (of length two) (cf. \cite{G,F, FP})
\begin{equation}
\begin{array}{l}
M_2(h)=-\int \omega\omega'.
\end{array}
\end{equation}
From previous calculation the only terms which contribute effectively are:
\begin{equation}
\begin{array}{l}
M_2(h)=-\int \omega_2\omega_5'-\int \omega_5\omega_2'.
\end{array}
\end{equation}
The main result we show here is that such an iterative integral can be computed by residues. For this purpose, we have first to compute $\omega_2, \omega_5$ and 
$\omega_2', \omega_5'$ in the coordinates $(z,h)$ and to determine their partial fraction decompositions.
\vskip 1pt
We recall an important formula (particular case of the shuffle formula see \cite{G}) for any couple of one-forms $\omega_0, \omega_1$:
\begin{equation}
\begin{array}{l}
\int \omega_0\omega_1+\int \omega_1\omega_0=\int \omega_0. \int \omega_1.
\end{array}
\end{equation}
In particular this yields that if $\int \omega_0=0$ or $\int \omega_1=0$, then
\begin{equation}
\begin{array}{l}
\int \omega_0\omega_1=-\int \omega_1\omega_0.
\end{array}
\end{equation}
\subsection{Computation of  $\omega_2$  and its derivatives}
We note that:
\begin{equation}
\omega_2=\frac{1}{2}hd(\ln(2y-1))=\frac{hdy}{2y-1}=\frac{x^2+y^2}{(2y-1)^2}dy,
\end{equation}
and thus we get:
\begin{equation}
\label{eq53}
\begin{array}{l}
\omega'_2=\frac{dy}{2y-1}=\\
\frac{1}{2{\rm i}z^2(2y-1)}(z^2+R^2)dz+\frac{[z-\frac{a+b}{2}]}{z(2y-1)}dh.
\end{array}
\end{equation}
If we change coordinates $(x,y)$ into $(z,h)$, we obtain:
\begin{equation}
\omega_2=h[\frac{z^2-ab}{2z(z-a)(z-b)}dz+{\rm i}\frac{(z-\frac{a+b}{2})}{(z-a)(z-b)}dh],
\end{equation}
and thus:
\begin{equation}
\omega_2=F_2(z,h)dz+\phi_2(z,h)dh,
\end{equation}
with
\begin{equation}
\begin{array}{l}
F_2(z,h)=\frac{h}{2}[-\frac{1}{z}+\frac{1}{z-a}+\frac{1}{z-b}]\\
\Phi_2(z,h)=\frac{{\rm i}h}{2}[\frac{1}{z-a}+\frac{1}{z-b}].
\end{array}
\end{equation}
We see that:
\begin{equation}
\begin{array}{l}
\int \omega_2=\int_{H=h} F_2 dz=0.
\end{array}
\end{equation}
The shuffle formula implies for instance:
\begin{equation}
\begin{array}{l}
\int \omega_2\omega_2=-\int \omega_2\omega_2=0.
\end{array}
\end{equation}
\subsection{Computation of $\omega_5$ and its derivatives in the coordinates $(z,h)$}
We change coordinates $(x,y)$ into $(z,H)$, this displays:
\begin{equation}
\begin{array}{l}
x=\frac{1}{2}\frac{(z^2+R^2)}{z}, y=\frac{-{\rm i}}{2}\frac{z^2+2{\rm i}hz-R^2}{z}\\
dx=\frac{(z^2-R^2)}{2z^2}dz+\frac{1}{2}\frac{2h-1}{z}dh, dy=\frac{1}{2{\rm i}z^2}(z^2+R^2)dz+\frac{[z-\frac{a+b}{2}]}{z}dh.
\end{array}
\end{equation}
We focus on:
\begin{equation}
\begin{array}{l}
\omega_5=\frac{2xydy-(x^2-y^2)dx}{(2y-1)^2}.
\end{array}
\end{equation}
We find:
\begin{equation}
\begin{array}{l}
\frac{2xy}{(2y-1)^2}dy=\frac{1}{4}\frac{(z^2-ab)^2[(z-a)(z-b)+{\rm i}z]}{[z(z-a)(z-b)]^2}dz\\
-\frac{1}{2{\rm i}}\frac{(z^2-ab)(z-\frac{a+b}{2})[(z-a)(z-b)+{\rm i}z]}{z[(z-a)(z-b)]^2}dh,
\end{array}
\end{equation}
\begin{equation}
\begin{array}{l}
\frac{-x^2+y^2}{(2y-1)^2}dx=\frac{1}{8}\frac{((z^2-ab)^2+[(z-a)(z-b)+{\rm i}z]^2)(z^2+ab)}{[z(z-a)(z-b)]^2}dz\\
+\frac{2h-1}{8}\frac{(z^2-ab)^2+[(z-a)(z-b)+{\rm i}z]^2}{z[(z-a)(z-b)]^2}dh.
\end{array}
\end{equation}
We compute now the partial fraction decomposition. The first term of the component factorizing $dz$ yields:
\begin{equation}
\begin{array}{l}
\frac{1}{4}\frac{(z^2-ab)^2}{[z^2(z-a)(z-b)]}+\frac{{\rm i}}{4}\frac{(z^2-ab)^2}{z[(z-a)(z-b)]^2}=\\
\frac{1}{4}\frac{(z^2-ab)^2}{[z(z-a)(z-b)]}[\frac{1}{z}-\frac{1}{z-a}+\frac{1}{z-b}]=\\
\frac{1}{4}[\frac{ab}{z}+\frac{a(a-b)}{z-a}-\frac{b(a-b)}{z-b}+a+b+z][\frac{1}{z}-\frac{1}{z-a}+\frac{1}{z-b}]=\\
\frac{1}{4}[1+2b[\frac{1}{z}-\frac{1}{z-a}+\frac{1}{z-b}]+\frac{ab}{z^2}-\frac{a(a-b)}{(z-a)^2}-\frac{b(a-b)}{(z-b)^2}].
\end{array}
\end{equation}
Then we consider:
\begin{equation}
\begin{array}{l}
\frac{1}{8}\frac{(z^2-ab)^2)(z^2+ab)}{[z(z-a)(z-b)]^2}=\\
\frac{1}{8}[1+2(a+b)\frac{1}{z}+\frac{2(a^2+b^2)}{a-b}[\frac{1}{z-a}-\frac{1}{z-b}]+\frac{ab}{z^2}+\frac{a(a+b)}{(z-a)^2}+\frac{b(a+b)}{(z-b)^2}].
\end{array}
\end{equation}
and
\begin{equation}
\begin{array}{l}
\frac{1}{8}\frac{[(z-a)(z-b)+{\rm i}z]^2(z^2+ab)}{[z(z-a)(z-b)]^2}=\\
\frac{1}{8}[1+\frac{ab}{z^2}]+\frac{{\rm i}}{4}\frac{z^2+ab}{[z(z-a)(z-b)]}-\frac{1}{8}\frac{z^2+ab}{[(z-a)(z-b)]^2}
\end{array}
\end{equation}
which gives for its partial fraction decomposition:
\begin{equation*}
\begin{array}{l}
\frac{1}{8}[1+2{\rm i}[\frac{1}{z}+\frac{a+b}{a-b}\frac{1}{z-a}-\frac{a+b}{a-b}\frac{1}{z-b}]+\\
\frac{4ab}{(a-b)^3}(\frac{1}{z-a}-\frac{1}{z-b})+\frac{ab}{z^2}\\
-\frac{a^2+ab}{(a-b)^2}\frac{1}{(z-a)^2}-\frac{b^2+ab}{(a-b)^2}\frac{1}{(z-b)^2}].
\end{array}
\end{equation*}
We define the two rational functions:
\begin{equation}
\label{eq66}
\omega_5=F_5dz+\Phi_5dh,
\end{equation}
and so the previous decomposition gives:
\begin{equation*}
\begin{array}{l}
F_5=\frac{1}{2}-{\rm i}(h-1)[\frac{1}{z}-\frac{1}{z-a}+\frac{1}{z-b}]\\
-\frac{h(h-1)}{2}[\frac{1}{z^2}+\frac{1}{(z-a)^2}]-\frac{(h-1)^2}{2}\frac{1}{(z-b)^2}.
\end{array}
\end{equation*}
From this, it is possible to compute by derivation:
\begin{align*}
F_{5h}'=&-{\rm i}[\frac{1}{z}-\frac{1}{z-a}+\frac{1}{z-b}]\\
&-(h-\frac{1}{2})\frac{1}{z^2}-\frac{1}{2}\frac{1}{(z-a)^2}-2(h-1)\frac{1}{(z-b)^2}\\
&+{\rm i}h(h-1)\frac{1}{(z-a)^3}+{\rm i}(h-1)^2\frac{1}{(z-b)^3}.
\end{align*}
The first contribution to $\Phi_5$ is
\begin{equation}
\frac{{\rm i}}{2}[\frac{(z^2+R^2)[(z+{\rm i}h)(z+{\rm i}(h-1))+{\rm i}z)](z+\frac{1}{2}{\rm i}(2h-1))}{z[(z+{\rm i}h)(z+{\rm i}(h-1))]^2}],
\end{equation}
and its partial fraction decomposition is:
\begin{align*}
\frac{1}{2}[{\rm i}+\frac{2h-1}{2}\frac{1}{z}-\frac{2h-1}{2}\frac{1}{z+{\rm i}h}+\frac{2h-3}{2}\frac{1}{z+{\rm i}(h-1)}\\
+\frac{{\rm i}h}{2}\frac{1}{(z+{\rm i}h)^2}+\frac{{\rm i}(h-1)}{2}\frac{1}{(z+{\rm i}(h-1))^2}].
\end{align*}
The second piece is:
\begin{equation}
\frac{(2h-1)}{8}[\frac{[(z^2+R^2)^2+[(z+{\rm i}h)(z+{\rm i}(h-1))+{\rm i}z]^2]}{z[(z+{\rm i}h)(z+{\rm i}(h-1)]^2}],
\end{equation}
and its partial fraction decomposition gives:
\begin{equation}
\frac{(2h-1)}{4}[\frac{1}{z}-\frac{1}{z+{\rm i}h}+\frac{1}{z+{\rm i}(h-1)}-\frac{{\rm i}h}{(z+{\rm i}h)^2}-\frac{{\rm i}(h-1)}{(z+{\rm i}(h-1))^2}].
\end{equation}
All together, the sum of the two pieces is:
\begin{equation}
\begin{array}{l}
\Phi_5=\frac{1}{2}[{\rm i}+\frac{2h-1}{z}-\frac{2h-1}{z+{\rm i}h}+\frac{2(h-1)}{z+{\rm i}(h-1)}\\
-\frac{{\rm i}h(h-1)}{(z+{\rm i}h)^2}-\frac{{\rm i}(h-1)^2}{(z+{\rm i}(h-1))^2}].
\end{array}
\end{equation}
\subsection{Computation of $ -\int \omega_5\omega'_2$ }
We begin by the observation that:
\begin{equation}
\omega'_2=\frac{1}{2}[-\frac{1}{z}+\frac{1}{z-a}+\frac{1}{z-b}]dz +...(dH),
\end{equation}
and so:
\begin{equation}
\int_{H=h} \omega'_2=0,
\end{equation}
hence we can apply the shuffle formula and obtain:
\begin{equation}
-\int \omega_5\omega'_2=\int \omega'_2\omega_5.
\end{equation}
This displays:
\begin{equation}
\begin{array}{l}
\int \omega'_2\omega_5=\int \frac{1}{2}[-\frac{1}{z}+\frac{1}{z-a}+\frac{1}{z-b}]dz\omega_5=\\
-\frac{{\rm i}(h-1)}{2}\int [-\frac{1}{z}+\frac{1}{z-a}+\frac{1}{z-b}]dz[\frac{1}{z}-\frac{1}{z-a}+\frac{1}{z-b}]dz\\
+\frac{1}{2} \int_{h=h}[-\frac{1}{z}+\frac{1}{z-a}+\frac{1}{z-b}]\{\frac{1}{2}z+h(h-1)[\frac{1}{z}+\frac{1}{z-a}]+\frac{(h-1)^2}{2}\frac{1}{z-b}\}dz.
\end{array}
\end{equation}
Note that the second expression can be readily computed by residue.
The first component breaks into four pieces that we compute by the shuffle formula:
\begin{equation}
\int [-\frac{1}{z}+\frac{1}{z-a}]dz[\frac{1}{z}-\frac{1}{z-a}]=0,
\end{equation}
\begin{equation}
\int [\frac{1}{z-b}]dz[\frac{1}{z}-\frac{1}{z-a}]dz=\int [-\frac{1}{z}+\frac{1}{z-a}]dz[\frac{1}{z-b}]dz,
\end{equation}
\begin{equation}
\int [\frac{1}{z-b}]dz[\frac{1}{z-b}]dz=0.
\end{equation}
This gives the contribution:
\begin{equation}
\begin{array}{l}
-{\rm i}(h-1)[\int_{H=h} (-\frac{1}{z}+\frac{1}{z-a}){\rm Log}(z-b) dz=\\
-{\rm i}(h-1)(2{\pi}{\rm i})[-{\rm Log}(-b)+{\rm Log}(a-b)]=\\
2{\pi}(h-1)[-{\rm Log}(-{\rm i}(1-h))+{\rm Log}(-{\rm i})]=
-2{\pi}(h-1){\rm Log}(1-h).
\end{array}
\end{equation}
The last component contributes to the sum of residues:
\begin{equation}
\frac{\pi}{2}h-\frac{{\pi}h(h-1)}{2}+\frac{{\pi}h}{2}-\frac{\pi}{2}(h-1)^2(\frac{1}{h-1}+1),
\end{equation}
and all together this holds:
\begin{equation}
-\int \omega_5\omega'_2=\pi(2h-h^2)-2\pi(h-1){\rm Log}(1-h).
\end{equation}
\subsection{Computation of $- \int \omega_2 \omega'_5$}
We begin with the expression of $\omega'_5$:
\begin{equation}
\omega'_5=(F'_{5h}-\Phi'_{5z})dz,
\end{equation}
and thus:
\begin{equation}
-\int \omega_2\omega'_5=-\int \omega_2F'_{5h} dz+\int F_{2}\Phi_5 dz.
\end{equation}
The second integral can be computed by residues and this yields:
\begin{equation}
\begin{array}{l}
\int_{H=h} \frac{h}{2}[-\frac{1}{z}+\frac{1}{z-a}+\frac{1}{z-b}].\\
\frac{1}{2}[{\rm i}+\frac{2h-1}{z}-\frac{2h-1}{z+{\rm i}h}+\frac{2(h-1)}{z+{\rm i}(h-1)}\\
-\frac{{\rm i}h(h-1)}{(z+{\rm i}h)^2}-\frac{{\rm i}(h-1)^2}{(z+{\rm i}(h-1))^2}].
\end{array}
\end{equation}
This term gives the contribution:
\begin{equation}
\begin{array}{l}
\frac{\pi}{2}\frac{h^2(2h-1)}{h-1}-{\pi}h(h-1)[\frac{1}{h-1}+1]+\frac{\pi}{2}h^2(h-1)-\frac{\pi}{2}h(h-1)^2+\frac{\pi}{2}h.
\end{array}
\end{equation}
We consider now:
\begin{equation}
\begin{array}{l}
-\int \omega_2 F'_{5h} dz={\rm i}\int \omega_2[\frac{1}{z}-\frac{1}{z-a}+\frac{1}{z-b}]dz\\
-\int \omega_2[(h-\frac{1}{2})\frac{1}{z}+\frac{1}{2}\frac{1}{z-a}+2(h-1)\frac{1}{z-b}\\
-{\rm i}\frac{h(h-1)}{2}\frac{1}{(z-a)^2}-{\rm i}\frac{(h-1)^2}{2}\frac{1}{(z-a)^2}]dz.
\end{array}
\end{equation}
The first term:
\begin{equation}
{\rm i}\int \omega_2[\frac{1}{z}-\frac{1}{z-a}+\frac{1}{z-b}]dz,
\end{equation}
can be computed using the shuffle formula as it was done in the previous paragraph and it yields:
\begin{equation}
\begin{array}{l}
\frac{{\rm i}h}{2}\int [-\frac{1}{z}+\frac{1}{z-a}+\frac{1}{z-b}]dz[\frac{1}{z}-\frac{1}{z-a}+\frac{1}{z-b}]dz=2{\pi}h{\rm Log}(1-h).
\end{array}
\end{equation}
The other terms can be computed by residue and their contribution is:
\begin{align*}
&-\pi\frac{h(h-\frac{1}{2})}{h-1}
+\frac{{\pi}h}{2}+2{\pi}h(h-1)[\frac{1}{h-1}
+1]-\frac{\pi}{2}h^2(h-1)\\
&-\frac{\pi}{2}h(h-1)^2[-1+\frac{1}{(h-1)^2}].
\end{align*}
The sum of the two contributions:
\begin{equation}
\begin{array}{l}
\frac{\pi}{2}\frac{h^2(2h-1)}{h-1}-{\pi}h(h-1)[\frac{1}{h-1}+1]+\frac{\pi}{2}h^2(h-1)-\frac{\pi}{2}h(h-1)^2+\frac{\pi}{2}h\\
-\pi\frac{h(h-\frac{1}{2})}{h-1}+\frac{{\pi}h}{2}+2{\pi}h(h-1)[\frac{1}{h-1}+1]\\
-\frac{\pi}{2}h^2(h-1)-\frac{\pi}{2}h(h-1)^2[-1+\frac{1}{(h-1)^2}]
\end{array}
\end{equation}
gives:
\begin{equation}
2{\pi}h^2.
\end{equation}

To conclude we have proved the:

\begin{theorem}
\label{thm2}
The value of $M_2(h)$ is:

\begin{align*}
M_2(h)&=\pi(2h-h^2)-2\pi(h-1){\rm Log}(1-h)+2\pi(h){\rm Log}(1-h)+2{\pi}h^2\\
&=2{\pi}[h+\frac{h^2}{2}+{\rm Log}(1-h)].
\end{align*}

\end{theorem}

\section{The second order Melnikov function : geometric computation}
\label{section3}
In this section we describe the geometric counterpart of the computations of the preceding section. By abuse of notations, denote by 
\begin{align}
\delta(h), \tilde\delta(h) \in \pi_1(\Gamma_h,*)
\end{align}
the  continuous families of simple loops whose homology classes (denoted with the same letter) were considered in section   \ref{firstmelnikov}, and such that
\begin{align*}
\delta(h) & = \{(x,y)\in \R^2 : \frac{x^2+y^2}{2y-1}=h\} , \; h<0 \\
\tilde \delta(h) & = \{(x,y)\in \R^2 : \frac{x^2+y^2}{2y-1}=h\} , \; h>1 .
\end{align*}

Assuming that $\int_{\delta(h)} \omega \equiv 0$, we define the second order Melnikov function $M_2$ associated to $(0,0)$ by the iterated integral
 $M_2(h)= \int_{\delta(h)} \omega \omega'$. If on the other hand  $\int_{\tilde \delta(h)} \omega \equiv 0$ then the second order Melnikov function associated to the center near $(0,1)$ is $\tilde M_2(h)= \int_{\tilde \delta(h)} \omega \omega'$. It is known $M_2$, $\tilde M_2$ satisfy a 
 linear differential equation of Fuchs-type \cite{G}. 
We are interested in the monodromy representation of $M_2$, $\tilde M_2$. For this purpose we need the orbit $\mathcal O_\delta$ of the closed loop $\delta(h)$, that is to say the set of free homotopy classes of loops, obtained from $\delta(h)$ by  "analytic continuation" with respect to the parameter $h$. The homology group  $H_1(\mathcal O_\delta)$ is just the set $\mathcal O_\delta$, but with a group structure, given by composition of loops.
There are many ways to compose two closed loops, but we consider two different compositions as representing the same element of $H_1(\mathcal O_\delta)$, which is achieved by taking a quotient with respect to the commutators $ [\mathcal O, \pi_1(\Gamma_{h_0},*)]$.

The main geometric fact about $M_2$ is that the map
$$
\delta \mapsto \int_{\delta} \omega \omega'
$$
is linear on $H_1(\mathcal O_\delta)$, that is to say
$$
\int_{\delta_1\circ \delta_2} \omega \omega' = \int_{\delta_1} \omega \omega' + \int_{\delta_2} \omega \omega' .
$$
Therefore the monodromy of $M_2$ is represented on $H_1(\mathcal O_\delta)$ which will be enough to deduce the explicit form of $M_2$, and therefore another proof of  Theorem \ref{thm2}. 

In section \ref{homology1} we compute $H_1(\mathcal O_\delta)$ and  $H_1(\mathcal O_{\tilde \delta})$.
This computation is independent from the rest of the paper. The computation of $M_2$, $\tilde M_2$ is carried out in section
\ref{section4}.

\subsection{The homology $H_1(\mathcal O_\delta)$  of the orbit $\mathcal O_\delta$ of the closed loop $\delta$}
\label{homology1}
The  first return map of $X_\lambda$, $\lambda\sim 0$,  constructed along a closed loop 
$$\delta=\delta(h)  \subset \{H(x,y)=h\}$$ of $X_0$ in a complex domain defines a germ of analytic automorphism $\C,0 \to \C,0$.
The dominant term of the return map with respect to parameters $\lambda$, the so called Bifurcation function, is an iterated path integral along $\delta(h)$, which depends on the free homotopy class of the closed orbit $\delta(h)$ in the leaf $\{H(x,y)=h\}$. The variation of $\delta(h)$ with respect to $h$ defines an orbit $\mathcal O_\delta$.
The monodromy representation of this Bifurcation function is then constructed on the so called "homology $H_1(\mathcal O_\delta)$ of the orbit of $\delta$" which is the purpose of the section. 

More precisely, let $\delta(h) \in \pi_1(\Gamma_h)$ be a free homotopy class of loops, depending continuously in the parameter $h$.
The fundamental group $\pi_1( \{h\in \C : h \neq 0, 1 \} ,h_0)$
of 
the set of regular values $$\{h\in \C : h \neq 0, 1 \}$$
 acts on the permutation group $Perm(\pi_1(\Gamma_{h_0}))$  of homotopy classes of closed loops on $\pi_1(\Gamma_{h_0})$
 and let $\mathcal O = \mathcal O_\delta \subset \pi_1(\Gamma_{h_0},*)$ be the smallest normal subgroup containing the orbit of $\delta(h_0)$ under this action. The homology of the orbit $\mathcal O$ is 
 $$
 H_1^\delta (\Gamma_{h_0}, \Z) = H_1(\mathcal O) \stackrel{def}{=} \mathcal O \slash (\mathcal O, \pi_1(\Gamma_{h_0},*) 
 $$
 where  $(\mathcal O, \pi_1(\Gamma_{h_0},*)$ is the commutator subgroup.
For details see  \cite{G,GI05}, where it was first defined. In what follows we use the notation $H_1(\mathcal O)$ introduced in \cite{mnop05} and we call it "homology of the orbit". 
The importance of  $H_1(\mathcal O)$ lies in the fact that the monodromy representation of the bifurcation function, in particular of the second Melnikov function $M_2$, is a sub-representation of 
$$
\pi_1(\{h\in \C : h \neq 0, 1 \},*) \to Aut (H_1(\mathcal O)) 
$$
this will be used in the next section.

The purpose of the present section is to compute the homologies
$$H_1(\mathcal O_\delta), H_1(\mathcal O_{\tilde \delta})$$ where 
$$
\delta(h), \tilde \delta(h) \subset \Gamma_h = \{ (x,y)\in \C^2 : \frac{x^2+y^2}{ 2y-1} = h \}
$$
are the continuous family of closed loops, vanishing at $(0,0)$ and $(0,1)$ when $h$ tends to $0$ or $1$, respectively. When $h\in \R$, the two families of closed loops 
$\delta(h), \tilde \delta(h)$ form the nests of periodic orbits shown on fig.\ref{fig1n}.  The Riemann sphere $\Gamma_h$ has two punctures over $y=1/2$, and denote the two simple loops making one turn around each of these punctures by $\alpha(h), \beta(h)$. We may suppose that  $\alpha(h), \beta(h)$ have a common starting point, so we can define the commutator
 $ (\alpha,\beta) =  \alpha^{-1}\beta^{-1} \alpha\beta$. The choice of this starting point will be irrelevant to the final result, which can be formulated as follows
\begin{theorem}
\label{homology}
The homology $H_1(\mathcal O_\delta)$ is a free $\Z$-module with two generators $\delta$ and $ (\alpha,\beta) =  \alpha^{-1}\beta^{-1} \alpha\beta$.
Similarly, the homology $H_1(\mathcal O_{\tilde \delta})$ is a free $\Z$-module with two generators $\tilde \delta$ and
$ (\alpha,\beta)$ .
\end{theorem}
The proof will be given later in this section.
Note that only the free homotopy classes of $\alpha, \beta, \delta, \tilde \delta $, specified on fig. \ref{figzf}, are relevant to the above statement.
The proof of Theorem \ref{homology}  is based on the algebraic Lemma \ref{algebraic} which we discuss first. For generalities on free groups see e.g. \cite{hall59}

Let $G$ be the free group generated by three letters $\alpha, \beta ,\gamma$. We consider the normal subgroup $H\subset G$ generated by the words
$$
\beta\gamma, (\alpha,\beta)
$$
where $(\alpha,\beta)= \alpha^{-1}\beta^{-1} \alpha \beta $ is the commutator of $\alpha,\beta$.

For  words  $x, w\in G$ we denote $x^w= w^{-1}x w$. 
For arbitrary words $w_1,w_2,w_3,w_4$, let $\tilde H\subset G$ be the normal subgroup of $G$ generated by the words
$$
\beta^{w_1}\gamma^{w_1}, (\alpha^{w_3},\beta^{w_4}) .
$$
\begin{lemma}
\label{algebraic}
With the above notations $H = \tilde H$. 
\end{lemma}
\proof
Let $(H,G)$ be the commutator subgroup of $H$, generated by commutators $(h,g)=h^{-1}g^{-1}hg$ where $h\in H, g\in G$.
For $w_1,w_2 \in G$
we have
$$
(\alpha^{w_1},\beta^{w_2})= (\alpha, \beta^{w_2w_1^{-1}})^{w_1}
$$
so we may suppose that the word $w_1$ is void, and consider the commutator $(\alpha,\beta^{w})$ for some $w\in G$. We have
\begin{align*}
(\alpha,\beta^{\beta})&= (\alpha,\beta) \\
(\alpha,\beta^{\alpha})&= (\alpha,\beta)^\alpha =  (\alpha,\beta) \mod (H,G) \\
(\alpha,\beta^{\gamma})&=(\alpha, \gamma^{-1}\beta\gamma) = (\alpha, \gamma^{-1}\gamma \beta)  \mod (H,G) \\
&=  (\alpha,  \beta)  \mod (H,G)
\end{align*}
It follows by induction, that for every word $w\in G$ holds
\begin{align}
\label{alphabeta}
(\alpha,\beta^{w}) =    (\alpha,\beta) \mod (H,G) 
\end{align}
and therefore $(\alpha,\beta^{w}) \in H$.
Similarly,$\beta^{w_1}\gamma^{w_2} =( \beta\gamma^{w_2w_1^{-1}})^{w_1}$ and we may assume as above
that $w_1=1$. We have
\begin{align}
\beta\gamma^{\beta}&= (\beta \gamma)^\beta = \beta\gamma \mod (H,G)  \\
\label{monbg}
\beta \gamma ^\alpha &=\alpha^{-1} \cdot \alpha \beta \alpha^{-1}\beta^{-1} \cdot \beta\gamma \cdot \alpha
= (\alpha,\beta) + \beta\gamma   \mod (H,G) \\
\beta\gamma^{\gamma}&= \beta\gamma .
\end{align}
It follows by induction that for every word $w\in G$ holds
\begin{align}
\label{betagamma}
\beta \gamma^w = \beta \gamma + k (\alpha,\beta)  \mod (H,G)  
\end{align}
where $k$ is an appropriate integer, and therefore $\beta \gamma^w \in H$.

Thus $\tilde H \subset H$ and in a similar way one shows that $H \subset \tilde H $ which implies $H = \tilde  H$.
\endproof
\begin{corollary}
According to Lemma \ref{algebraic} the free Abelian factor group $H/(H,G)$ generated by $\beta\gamma$ and the commutator $(\alpha,\beta)$ depends in fact only on the conjugacy classes of the letters $\alpha, \beta, \gamma$. 
\end{corollary}

We shall apply Lemma \ref{algebraic} in the following geometric situation.

Using the notations of section \ref{firstmelnikov}, let $\Gamma_h$ is the four-punctured Riemann sphere identified to the complex $z$-plane $\C$, with coordinate
$z= x+i(y-h)$ and punctures at 
$
a= -ih, b =  -i(h-1), c = 0 .
$
It is seen that $a,b\in \C$ are arbitrary constants subject to the relation $b-a=i$, see fig. \ref{figzf}.

\begin{figure}
\includegraphics[width=8cm]{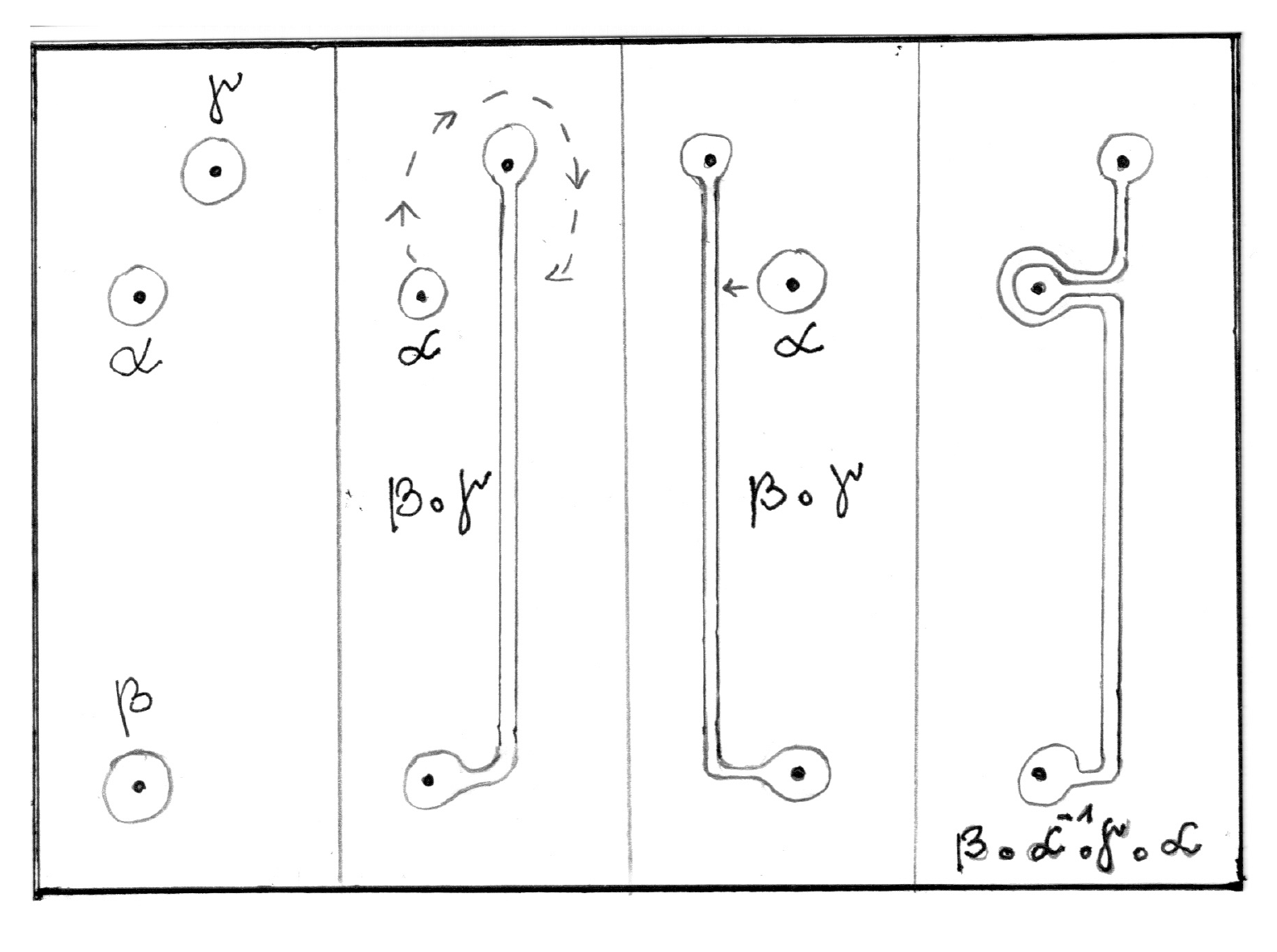}
\caption{The monodromy   of the closed loop $\delta(h)=\beta(h)\circ \gamma(h)$ when $h$ makes one turn around $h=1$.}
\label{monodromy}
\end{figure}

Let 
$$\alpha, \beta,\gamma, \delta = \gamma\alpha, \; \tilde \delta = \beta\gamma  \in \pi_1(\Gamma_h, *)$$
be such, that their corresponding free homotopy classes are represented by the closed loops denoted with the same letter on fig. \ref{figzf}.
As in Lemma \ref{algebraic}, let $G$ be the group generated by $\alpha, \beta,\gamma$, $H$ its normal subgroup generated by 
$\tilde \delta, (\alpha,\beta)$. 
\proof[Proof of Theorem \ref{homology}]
We shall prove that $\mathcal O_{\tilde \delta} = H$. 
When $h\in (0,1)$ is close to $1$, and then makes one turn around $h=1$ along the path
$$
h \mapsto 1+e^{i\varphi} (h-1),  \; \varphi \in [0,2 \pi ]
$$
the resulting free homotopy class $ \tilde \delta(1+e^{2 \pi i} (h-1))$ equals $ \tilde \delta(h)$, so its variation is trivial.  When $h\in (0,1)$ is close to $0$, and then makes one turn around $h=0$ along the path
$$
h \mapsto e^{i\varphi} h,  \; \varphi \in [0,2 \pi ]
$$
we find, according to fig.\ref{monodromy} that $\tilde \delta = \beta\gamma$ undergoes the following monodromy transformation
$$
\tilde \delta = \beta\gamma \mapsto \beta \alpha^{-1} \gamma \alpha 
$$
and according to (\ref{betagamma})
$$
\beta \alpha^{-1} \gamma \alpha  =\alpha^{-1} \cdot \alpha \beta \alpha^{-1}\beta^{-1} \cdot \beta\gamma \cdot \alpha
= (\alpha,\beta) + \beta\gamma   \mod (H,G) .
$$
It remains to compute the monodromy of the commutator $(\alpha,\beta) $. By analogy to fig.\ref{monodromy} we find  that 
$(\alpha,\beta)$ undergoes the transformation
$$
(\alpha,\beta) \mapsto (\alpha, \gamma^{-1} \beta \gamma)
$$
and according to (\ref{alphabeta})
$$
(\alpha, \gamma^{-1}\beta\gamma) =  (\alpha,  \beta)  \mod (H,G) .
$$
This completes the proof that $\mathcal O_{\tilde \delta} = H$ and hence 
$$
H_1(\mathcal O_{\tilde \delta}) = H/(H,G) .
$$
The computation of $H_1(\mathcal O_\delta)$ repeats the same arguments. \endproof

\subsection{The monodromy of the second Melnikov function }
\label{section4}
Recall that, under the condition that the first Melnikov function $ \int_{\delta(h)} \omega$ is identically zero, we have for the second
$M_2(h) = 4 \int_{\delta} \omega \omega'$ where $\delta = \beta \circ \gamma$ as it is shown on 
fig.\ref{monodromy}. Clearly, $M_2$ is analytic near $h=0$ and analytic on the domain
$$
\{h \in \C \} \setminus [1, \infty) .
$$
It allows  an analytic continuation along the universal covering of $\C~\setminus\{0,1\}$.
When $h$ makes one turn around $h=1$, according to fig.\ref{monodromy} and formula (\ref{monbg}) the monodromy of $\delta$ as an element of $H_1(\mathcal O_\delta) = \Z\delta + \Z(\alpha,\beta)$ is
$$
\delta \mapsto \delta + (\alpha,\beta) .
$$
On the other hand  $(\alpha(h),\beta(h)$ have no monodromy at all, see (\ref{alphabeta}). Therefore 
$$
 \int_{\delta(h)} \omega \omega' = P(h) + Q(h) \ln (1-h)
$$
where $P, Q$ are rational functions in $h$ with eventual pole at $h=1$, and moreover
$$
2\pi i Q(h) = \int_{(\alpha(h),\beta(h))} \omega \omega' .
$$
Of course, similar considerations are valid for the family of loops $\tilde \delta (h)$. Namely, under the condition that $\int_{\tilde \delta(h)} \omega = 0$ we have
$$
 \int_{\tilde \delta(h)} \omega \omega' = \tilde P(h) + \tilde Q(h) \ln (h)
$$
where $\tilde P, \tilde Q$ are rational functions in $h$ with eventual pole at $h=0$, and moreover
$$
2\pi i \tilde Q(h) = \int_{(\alpha(h),\beta(h))} \omega \omega' .
$$

The iterated integral along the commutator $(\alpha,\beta) = \alpha^{-1}\beta^{-1}\alpha \beta$ is however easily computed by standard properties of iterated integrals, e.g.
 \cite[Lemma A.2]{G}, so  we have
\begin{align}
\label{albe}
\int_{(\alpha,\beta)} \omega \omega' =
 \det \left(\begin{array}{cc}\int_{\alpha}\omega  & \int_{\alpha}\omega' \\ \int_\beta \omega & \int_{\beta}\omega' \end{array}\right) .
\end{align}
By (\ref{eq53}) and (\ref{eq66}) we have that along $\Gamma_h$
\begin{align*}
\omega_2 &= \frac{h}{2}(-\frac{1}{z}+\frac{1}{z-a}+\frac{1}{z-b}) dz 
\end{align*}
\begin{align*}
\omega_5 &=  [\frac{1}{2}-{\rm i}(h-1)(\frac{1}{z}-\frac{1}{z-a}+\frac{1}{z-b}) \\
&-\frac{h(h-1)}{2}(\frac{1}{z^2}+\frac{1}{(z-a)^2})-\frac{(h-1)^2}{2}\frac{1}{(z-b)^2}] dz .
\end{align*}
Recall that $\alpha, \beta$ a simple loops around $a=-ih$ and $b-i(h-1)$ respectively.
Residue calculus  implies that for $ \omega= \lambda_2\omega_2+\lambda_5\omega_5$ holds
\begin{align*}
\frac{1}{2\pi i} \int_\alpha \omega = \frac{h\lambda_2}{2} + i(h-1)\lambda_5, 
\frac{1}{2\pi i} \int_\beta \omega = \frac{h\lambda_2}{2} -  i(h-1)\lambda_5
\end{align*}
and taking into account that 
\begin{align*}
\int_{\alpha(h)} \omega'= \frac{d}{dh}\int_{\alpha(h)}\omega, 
\int_{\beta(h)} \omega'= \frac{d}{dh}\int_{\beta(h)}\omega
\end{align*}
we conclude that
\begin{align}
\label{96}
\int_{(\alpha,\beta)} \omega \omega' = -4\pi^2 i \lambda_2\lambda_5 .
\end{align}
Therefore
\begin{align*}
 \int_{\delta(h)} \omega \omega' = P(h) - 2 \pi \lambda_2 \lambda_5 \ln(1-h) .
\end{align*}
As expected the coefficient of $\ln (1-h)$ is quadratic in $\lambda_i$ and is therefore a generator of $\mathcal B_1$, that is to say $\lambda_2\lambda_5$.
In particular it should not contain $\lambda_4$ so the formula (\ref{96}) is also valid for
$$
\omega = \lambda_2 \omega_2 + \lambda_4 \omega_4 + \lambda_5 \omega_5 .
$$
Similarly, the condition $\tilde M_1=0$ implies $\lambda_1+\lambda_3=\lambda_5=0$ and hence for 
$$
\omega =  \lambda_2\omega_2 + \lambda_3 (\omega_3-\omega_1) +\lambda_4 \omega_4
$$
\begin{align}
\label{97}
\int_{(\alpha,\beta)} \omega \omega' = -4\pi^2 i \lambda_3\lambda_4 
\end{align}
which on its turn implies
\begin{align*}
 \int_{\tilde \delta(h)} \omega \omega' = \tilde P(h) - 2 \pi \lambda_3 \lambda_4 \ln(h) .
\end{align*}

To compute $P(h)$ we use asymptotic analysis. It can be verified that the meromorphic function  $P(h)$ has no poles on the finite plane, and grows at infinity no faster than $h^2$ (we skip the proof). From this already follows that
$P(h)$ is a polynomial of degree at most two. To find its  coefficients we use the following
\begin{proposition}
\label{expansion}
The second Melnikov function $M_2= 4  \int_{\delta(h)} \omega \omega' $ has a zero at $h=0$ of multiplicity at least three.
\end{proposition}
Indeed, assuming the Proposition, by the expansion $$\ln(1-h) = -h- \frac{h^2}{2} +\dots$$ we get $P(h)$, and hence the main result of this section
\begin{align}
\boxed{
 \int_{\delta(h)} \omega \omega' = 2\pi \lambda_2 \lambda_5  (h + \frac{h^2}{2}+ \ln(1-h)) }
\end{align}
which agrees, as expected,  with  Theorem \ref{thm2}.

\proof[Proof of Proposition \ref{expansion}]
A local analytic change of the variables in a neighborhood of $(0,0)$
\begin{align}
\label{morse}
\frac{x}{\sqrt{1-2y}} \mapsto x,\; \frac{y}{\sqrt{1-2y}} \mapsto y
\end{align}
brings the Hamiltonian $H= \frac{x^2+y^2}{2y-1}$ to the form $H(x,y) = -x^2-y^2$. To compute 
$$
\int_{\{H=h\}} \omega \omega'
$$
we may suppose that $\lambda_1=\lambda_3=\lambda_4=0$. The differential
\begin{align*}
\omega =  \lambda_2  \frac{(x^2  +y^2 )dy}{(2y-1)^2}
 + \lambda_5   \frac{  2xy dy - (x^2-y^2)dx }{(2y-1)^2} .
\end{align*}
in the new coordinates (\ref{morse}) takes the form
\begin{align*}
\omega =  -\lambda_2  (x^2+y^2) dy
 + \lambda_5   d ( xy^2 - \frac{x^3}{3}) + O(4)
\end{align*}
where the $O(4)$ replaces some analytic differential one-form $Pdx+Qdy$, where $P,Q,$ vanish at the origin$(0,0)$ of multiplicity at least three.
We have therefore
\begin{align*}
\omega =&  -\lambda_5  d(H y) + \lambda_5   d ( xy^2 - \frac{x^3}{3}) + \lambda _2y dH + \dots \\
\omega' =& \lambda _2 dy + \dots
\end{align*}
and

\begin{align*}
 \int_{\{H=h\}} \omega \omega' = &  \lambda_2 \int_{\{H=h\}}  y \omega  + O(5) \\
=& \lambda_2 \lambda_5 \int_{\{H=h\}} y \; d ( xy^2 - \frac{x^3}{3}) +O(5)\\
=& \iint_{ \{x^2+y^2 < -h\}} (y^2 - x^2)dy\,dx + O(5) = O(5).
\end{align*}
As $\deg H = 2$, then homogeneity considerations show that $O(5)=O(h^3)$ so $M_2(h)=  O(h^3)$
\endproof

\section{Blow up of a direct product of ideals}
\label{sectionblowup}

Let $\mathbb C
\{{ \lambda }\}$ be the ring of convergent power series at $\lambda=0$, where
$(\lambda_1,...\lambda_n)\in  \C^n$, and 
$$\mathcal B = (v_1,\dots,v_N) \subset \mathbb C
\{{ \lambda }\}$$
 be an ideal with zero set
 $$
 Z(\mathcal B)= \{ \lambda \in (\C^n,0) : v_1(\lambda)=v_2(\lambda)=\dots=v_N(\lambda)=0 \}
 $$ 
The blowup $\Gamma_\mathcal B \subset (\C^n,0) \times \mathbb P^{N-1}$
 of $(\C^n,0)$
with center $\mathcal B$ is the analytic closure
 of the graph of the map
\begin{align*}
C^n \setminus  Z(\mathcal B) &\to \mathbb P^{N-1}\\
{\lambda }  &\mapsto
[ v_1({\lambda }):\dots : v_N({\lambda })]
\end{align*}
with projection on the first factor 
$$\pi_\mathcal B  : \Gamma_\mathcal B \subset (\C^n,0) \times \mathbb P^{N-1}  \to (\C^n,0) .$$
Here $ [ v_1({\lambda }):\dots : v_N({\lambda })] $
is the projectivization of $ (
v_1({\lambda }),\dots , v_N({\lambda }))$.
The exceptional divisor
$$
E_\mathcal B=\pi^{-1}(0) \subset  \mathbb P^{N-1}
$$
is therefore a well defined closed algebraic set. The importance of $E_\mathcal B$ lies in the fact that
it is in bijective correspondence with the projectivized set of bifurcation (or Melnikov)  functions, computed in the preceding sections, see \cite[Corollary 2]{FGX}.

Suppose that $\mathcal B_1, \mathcal B_2 \subset \mathbb C \{{ \lambda }\} $ be two ideals 
\begin{align*}
 \mathcal B_1 & = (v_1^1,\dots,v_{N_1}^1) \\
 \mathcal B_2 & = (v_1^2,\dots,v_{N_1}^2) 
\end{align*}
and consider
the direct product
$$
\mathcal B = \mathcal B_1 \times \mathcal B_2 \subset \mathbb C \{{ \lambda }\} \times \mathbb C \{{ \lambda }\}.
$$
We note that $\mathcal B $ is also an ideal and consider the corresponding blowup
$$
\Gamma_\mathcal B \subset (\C^n,0) \times \mathbb P^{N_1-1} \times  \mathbb P^{N_2-1} 
$$
defined as the analytic closure of the graph of the map
\begin{align*}
C^n \setminus  Z(\mathcal B) &\to \mathbb P^{N_1-1} \times  \mathbb P^{N_2-1} \\
{\lambda }  &\mapsto
([ v_1^1({\lambda }):\dots : v_{N_1}^1({\lambda })], [ v_1^1({\lambda }):\dots : v_{N_2}^1({\lambda })] ) 
\end{align*}
with corresponding exceptional divisor
$$
E_{\mathcal B_1 \times \mathcal B_2 }=\pi^{-1}(0) \subset  \mathbb P^{N_1-1} \times  \mathbb  P^{N_2-1} .
$$
To the end of the present section we compute $E_{\mathcal B_1 \times \mathcal B_2 }$ in the case when
\begin{align*}
\mathcal B_1 &= < v_1^1(\lambda),v_2^1(\lambda),v_3^1(\lambda)> =<\lambda_1,\lambda_3, \lambda_2\lambda_5 >\\
\mathcal B_2 &=  < v_1^2(\lambda),v_2^2(\lambda),v_3^2(\lambda)> =<\lambda_1+\lambda_3 + \lambda_1 \lambda_2,\lambda_5, \lambda_3 \lambda_4> .
\end{align*}
It follows with same proof as \cite[Corollary 2]{FGX} that
\begin{proposition}
The projectivized set of pairs of Melnikov functions computed in the preceding section are in bijective correspondence with the points on the exceptional divisor $E_{\mathcal B_1 \times \mathcal B_2 }$.
\end{proposition}
The main result of the present section is
\begin{theorem}
\label{blowup}
The exceptional divisor  $$E_{\mathcal B_1 \times \mathcal B_2 } \subset \mathbb P^2 \times \mathbb P^2$$
has three irreducible components as follows
\begin{align}
\label{set1}
\{ ([c_1^1:c_2^1:c_3^1],[ c_1^2:c_2^2:c_3^2]): c_1^2=c_3^2=0 \}   \\
\label{set2}
\{ ([c_1^1:c_2^1:c_3^1],[ c_1^2:c_2^2:c_3^2])  : c_1^1+c_2^1=0, c_3^1=0 \}   \\
\label{set3}
\{ ([c_1^1:c_2^1:c_3^1],[ c_1^2:c_2^2:c_3^2]): c_3^1=c_3^2=0 \}  .
\end{align}
\end{theorem}
\proof
A point $(P_1,P_2) \in (\mathbb P^2,\mathbb P^2)$ belongs to $E_{\mathcal B_1 \times \mathcal B_2 }$ if and only if there is an arc
\begin{align}
\label{arc1}
\varepsilon \mapsto \lambda(\varepsilon) = (\lambda_1(\varepsilon), \dots, \lambda_6(\varepsilon) ), \lambda(0)=0
\end{align}
such that the vector
$$
([v_1^1(\lambda(\varepsilon)):v_2^1(\lambda(\varepsilon)):v_3^1(\lambda(\varepsilon))],[ v_1^2(\lambda(\varepsilon)):v_2^2(\lambda(\varepsilon)):v_3^2(\lambda(\varepsilon))])
$$
tends to the vector $(P_1,P_2)$ as $\varepsilon$ tends to $0$. It is easy to show now that 
the components (\ref{set1}),(\ref{set2}), (\ref{set3}) belong to $E_{\mathcal B_1 \times \mathcal B_2 }$. For instance, for  (\ref{set2}) we may consider the family of arcs
$$
\varepsilon \mapsto
\lambda( \varepsilon)=
(\varepsilon, \varepsilon^2, -\varepsilon + \lambda_3^0\varepsilon^2, -\lambda_4^0\varepsilon, \lambda_5^0\varepsilon^2)
$$
and then 
\begin{align*}
\lim_{\varepsilon \mapsto 0} &([v_1^1(\lambda(\varepsilon)):v_2^1(\lambda(\varepsilon)):v_3^1(\lambda(\varepsilon))],[ v_1^2(\lambda(\varepsilon)):v_2^2(\lambda(\varepsilon)):v_3^2(\lambda(\varepsilon))]) \\
& = ([1:-1:0],[\lambda_3^0:\lambda_5^0:\lambda_4^0])
\end{align*}
The other inclusion are also obvious.

Next, we consider an arbitrary arc (\ref{arc}) and we must show that 
$$
\lim_{\varepsilon \mapsto 0} ([v_1^1(\lambda(\varepsilon)):v_2^1(\lambda(\varepsilon)):v_3^1(\lambda(\varepsilon))],[ v_1^2(\lambda(\varepsilon)):v_2^2(\lambda(\varepsilon)):v_3^2(\lambda(\varepsilon))]) 
$$
belongs to one of (\ref{set1}),(\ref{set2}), (\ref{set3}). For this purpose we note that for fixed $\lambda_2, \lambda_4$, the generators $v_i^j$ of
$\mathcal B_1 , \mathcal B_2 $ are linear homogeneous in $\lambda_1, \lambda_3, \lambda_5$.
Therefore we shall consider separately each of the cases
$$
d_1 = \min_{i=1,3,5} d_i, d_3 = \min_{i=1,3,5} d_i, d_5 =  \min_{i=1,3,5} d_i .
$$
where
$$\lambda_1= O(\varepsilon^{d_1}), \lambda_3= O(\varepsilon^{d_3}), \lambda_5= O(\varepsilon^{d_5}) .$$
\begin{itemize}
\item The case $d_1 = \min_{i=1,3,5} d_i$
We put 
$$\lambda_1 = \lambda_1^0 \varepsilon^{d_1} + \dots,
\lambda_3 = \lambda_3^0 \varepsilon^{d_3} + \dots,
\lambda_5= \lambda_5^0 \varepsilon^{d_5} + \dots
$$
and hence
$$
\lim_{\varepsilon \mapsto 0} [\lambda_1:\lambda_3: \lambda_2\lambda_5 ]= [\lambda_1^0:\lambda_3^0:0] .
$$
If $\lambda_1+\lambda_3 = O(\varepsilon^{d_1})$ then 
\begin{align*}
\lim_{\varepsilon \mapsto 0} [\lambda_1+\lambda_3 + \lambda_1 \lambda_2:\lambda_5: \lambda_3 \lambda_4]=
[*,*,0] 
\end{align*}
and therefore the limit is in the set (\ref{set3}). If, however $\lambda_1+\lambda_3 = O(\varepsilon^{\tilde d_1})$ where $\tilde d_1 > d_1$, then
$$
\lim_{\varepsilon \mapsto 0} [\lambda_1:\lambda_3: \lambda_2\lambda_5 ]= [1:-1:0] .
$$
and the limit is in the set (\ref{set2}).

\item The case $d_3 = \min_{i=1,3,5} d_i$
We may suppose in addition that $d_3<d_1$ (otherwise we are in the preceding case). Then we check immediately that the limit is in the set (\ref{set3})
\item The case $d_5 = \min_{i=1,3,5} d_i$
We may suppose in addition that $d_5<d_1$ and $d_5<d_3$ (otherwise we are in one of the preceding two cases). Therefore 
\begin{align*}
\lim_{\varepsilon \mapsto 0} [\lambda_1+\lambda_3 + \lambda_1 \lambda_2:\lambda_5: \lambda_3 \lambda_4]=
[0:1:0] 
\end{align*}
and we are in the case (\ref{set1}).
\end{itemize}
This completes the proof of Theorem \ref{blowup}.
\endproof

\section{Distributions of limit cycles}
In this section we determine the possible distributions $(i,j)$ of limit cycles of small quadratic deformations (\ref{pert1}) of the quadratic vector field
(\ref{nona}) on the finite plane $\R^2$. This excludes the limit cycles, which bifurcate from "infinity".

\begin{definition} We  say that the germ of a family of vector fields $X_{a,b}$ 
\begin{eqnarray*} X_{a,b} :\left\{\begin{aligned}
\dot{x}=&   - y -x^2+y^2 + \sum_{0\leq i,j\leq 2 } a_{ij} x^iy^j ,\\
\dot{y}= & \;\;\;\;x - 2xy  -  \sum_{0\leq i,j\leq 2 }   b_{ij} x^iy^j 
\end{aligned}
\right.
\end{eqnarray*}
has an
admissible distribution $(i,j)$ of limit cycles, if there is a sequence $(a_k,b_k)_k$ in the parameter space $\{(a,b)\}$ such that for every sufficiently big $R\in \R$ the following holds true :
every vector field $X_{a_k,b_k}$ has exactly $i$  limit cycles surrounding the equilibrium point near $(0,0)$, exactly $j$  limit cycles surrounding the equilibrium point near $(0,1)$, and these limit cycles are  contained in the disc $\{(x,y) \in \R^2: \|(x,y\| <R \}$.
\end{definition}
The maximal value of $i$ is therefore the cyclicity  $Cycl(\Pi_1, X_{a,b})$ of the open period annulus containing $(0,0)$, the  maximal value of $j$ is the cyclicity  $Cycl(\Pi_2, X_{a,b})$ of the open period annulus containing $(0,1)$, and finally
$$
\max_{i,j} i+j = Cycl(\R^2, X_{a,b})
$$
Recall that the cyclicity  $Cycl(\Pi, X_{a,b})$ of an open set $\Pi\subset \R^2$ 
with respect to the germ of a family of vector fields $X_{a,b}$
is, roughly speaking, the maximal number of limit cycles which bifurcate from an arbitrary compact set $K\subset \Pi$ when $a,b \sim 0$. For a precise definition see e.g. \cite[Definition 3]{G01}. 

The main result of the paper is
\begin{theorem}
\label{main}
The distribution $(i,j)$ of limit cycles is admissible if and only if $i+j\leq 2$.
\end{theorem}
\proof
Without loss of generality we replace the germ of families  $X_{a,b}$ by $X_\lambda$, see (\ref{perturbed2}). The first return maps $\mathcal P_1, \mathcal P_2$ parameterized by the restriction $h=H(x,y)$ of the first integral on a cross-section to the annulus $\Pi_1$ or $\Pi_2$ can be divided 
in the corresponding ideals (\ref{bb1}) and (\ref{bb2}) as follows
\begin{align*}
\mathcal P_1^1(h;\lambda)(h) - h &=  v_1^1(\lambda)(M_1^1(h) + O(\lambda)) + v_2^1(\lambda)(M_2^1(h) 
+ O(\lambda)) \\ 
&+ v_3^1(\lambda)(M_3^1(h) + O(\lambda))\\
\mathcal P_1^2(h;\lambda)(h) - h &= v_1^2(\lambda)(M_1^2(h) + O(\lambda)) + v_2^2(\lambda)(M_2^2(h) + O(\lambda)) \\
&+v_3^2(\lambda)(M_3^2(h) + O(\lambda))
\end{align*}
where the Melnikov functions $M_i^j$ were computed in the preceding sections.
It follows, with same proof as    \cite[Theorem 1]{G08} that if $(i,j)$ is an admissible distribution of limit cycles for $X_\lambda$, then there exists a germ of analytic arc $$\varepsilon \mapsto \lambda(\varepsilon), \; \varepsilon \in (\R,0), \;\lambda(0)=0 $$ 
such that the one-parameter family of vector fields $X_{\lambda(\varepsilon)}$ allows a distribution $(i,j)$ of limit cycles, for  $\varepsilon$ close to $0$. For such an arc we obtain
\begin{align*}
\mathcal P_1^1(h;\lambda(\varepsilon))(h) - h &= \varepsilon^{k_1}( c_1^1 M_1^1(h)  + c_2^1 M_2^1(h) 
+   c_3^1 M_3^1(h) + O(\varepsilon))\\
\mathcal P_1^2(h;\lambda(\varepsilon))(h) - h &=  \varepsilon^{k_2}( c_1^2 M_1^2(h)  + c_2^2 M_2^2(h) 
+   c_3^2 M_3^2(h) + O(\varepsilon))
\end{align*}
Therefore to compute the distribution $(i,j)$ of limit cycles we have to compute the number of zeros $i$ and $j$
of each admissible pair of bifurcation functions
$$
c_1^1 M_1^1(h)  + c_2^1 M_2^1(h) 
+   c_3^1 M_3^1(h) , c_1^2 M_1^2(h)  + c_2^2 M_2^2(h) 
+   c_3^2 M_3^2(h) 
$$
According to section \ref{blowup} and 
(\ref{m1}), (\ref{m1t}),  the bifurcation function  associated to the first period annulus is co-linear to
$$
h[c_1^1(h-1) + c_2^1 h + c_3^1 M_2(h)]
$$
and  the bifurcation function associated to the second annulus is 
$$
(h-1) [c_1^2(h-1) -2 c_2^2  + c_3^2 \tilde M_2(h)] .
$$

According to \cite[Corollary 2]{FGX}
the admissible pairs of vectors 
$$
[c_1^1:c_2^1:c_3^1] , [c_1^2:c_2^2:c_3^2]  \in \mathbb P^2  
$$
are in one-to-one correspondance to the points on the exceptional divisor $$E_{\mathcal B_1 \times \mathcal B_2 } \subset \mathbb P^2 \times \mathbb P^2$$ described in Theorem \ref{blowup}. We consider each of the three irreducible components of $E_{\mathcal B_1 \times \mathcal B_2 } $ separately.

In the component (\ref{set1}) we have  
$c_1^2=c_3^2=0 $ so the bifurcation function associated to the second annulus $\Pi_2$ is co-linear to
$
h-1 .
$. Thus no limit cycles bifurcate from $\Pi_2$ and at most two limit cycles bifurcate from $\Pi_1$.

In the component (\ref{set2}) we have  
$c_1^1+c_2^1=0, c_3^1=0$ and hence the bifurcation function associated to the first period annulus $\Pi_1$ is co-linear to $h$.  Thus no limit cycles bifurcate from $\Pi_1$ and at most two limit cycles bifurcate from $\Pi_2$.

In the component (\ref{set3}) we have  
$ c_3^1=c_3^2=0$ and hence the bifurcation functions associated to the period annuli are co-linear to
\begin{align*}
h[c_1^1(h-1) + c_2^1 h], \; (h-1) [c_1^2(h-1) -2 c_2^2 ] .
\end{align*}
Therefore in each period annulus at most one limit cycle can bifurcate. This completes the proof.
\endproof

\section{Conclusion and Perspectives}

In this article, we have computed the double Bautin ideal associated to the bifurcations of the Lotka-Volterra double center with respect to arbitrary quadratic deformations. Our approach is based on the expression of the second-order bifurcation function in terms of iterated path integrals and on the shuffle formula. We also provide a geometric approach based on the ``homology of the orbit" description of the monodromy of the second-order bifurcation function of independent interest. Our main result allows to prove that the second-order bifurcation function is enough to compute the double cyclicity. Although we recall that the methods we have used do not allow to keep track of all the limit cycles which are born at the boundaries of the period annuli. This issue has been adressed in several other bifurcation settings \cite{CLY, DLZ, GI05, GI15}. This is certainly an interesting perspective for further researches. Another important perspective would be to try to extend the
outline of a general bifurcation theory of plane systems of infinite co-dimension that we have introduced here, in particular to other reversible quadratic double centers.

\bibliographystyle{plain}

\end{document}